\documentclass[12pt]{amsart}
\usepackage{graphics}
\usepackage{latexsym,bm}
\usepackage{amsmath,amsfonts,amssymb}

\newcommand{\BS}{\mathfrak S}
\newcommand{\Z}{\mathbb Z} \newcommand{\ui}{\underline i}
\newcommand{\ub}{\underline b} \newcommand{\uk}{\underline k}
\newcommand{\uj}{\underline j}
\newcommand{\T}{\mathcal T} \newcommand{\PP}{\mathcal P}
\newcommand{\CM}{\mathcal M}
\newcommand{\D}{\mathfrak D}
\newcommand{\Q}{\mathbb Q} \newcommand{\lam}{\lambda}
 \DeclareMathOperator{\End}{End}
\newcommand{\ft}{\mathfrak t} \DeclareMathOperator{\Ind}{Ind}
\DeclareMathOperator{\Ker}{Ker} \DeclareMathOperator{\res}{res}
\DeclareMathOperator{\Std}{Std} 
\DeclareMathOperator{\fa}{\bf a} \DeclareMathOperator{\fb}{\bf b}
 \DeclareMathOperator{\RS}{RowStd}
 \DeclareMathOperator{\Span}{Span}
\DeclareMathOperator{\GL}{GL}  \DeclareMathOperator{\fS}{S}
\DeclareMathOperator{\fT}{T} \DeclareMathOperator{\fI}{I}
 \DeclareMathOperator{\BD}{Bd}
\newcommand{\fs}{\mathfrak s} \newcommand{\fu}{\mathfrak u}

 \newcommand{\sig}{\sigma}

\newcommand{\bb}{\mathfrak{B}}

\newtheorem{prop}{Proposition}
\newtheorem{thm}{Theorem}\newtheorem{cor}{Corollary}
\newtheorem{lem}{Lemma}\newtheorem{dfn}{Definition}
\newtheorem{rem}{Remark}
\numberwithin{equation}{section} \numberwithin{prop}{section}
\numberwithin{thm}{section} \numberwithin{lem}{section}
\numberwithin{dfn}{section} \numberwithin{cor}{section}
\numberwithin{rem}{section}

\begin{document}
\title[Specht filtrations and tensor spaces]{Specht filtrations and tensor spaces for the Brauer algebra}
\author{Jun Hu}
\address{Department of Applied Mathematics,
Beijing Institute of Technology,
Beijing 100081 P.R. China}
\email{junhu303@yahoo.com.cn}

\thanks{Research supported by National Natural Science Foundation of
China (Project 10401005) and by the Program NCET}

\begin{abstract}
Let $m, n\in{\mathbb N}$. In this paper we study the right
permutation action of the symmetric group $\BS_{2n}$ on the set of
all the Brauer $n$-diagrams. A new basis for the free $\Z$-module
$\bb_n$ spanned by these Brauer $n$-diagrams is constructed, which
yields Specht filtrations for $\bb_n$. For any $2m$-dimensional
vector space $V$ over a field of arbitrary characteristic, we give
an explicit and characteristic free description of the annihilator
of the $n$-tensor space $V^{\otimes n}$ in the Brauer algebra
$\bb_n(-2m)$. In particular, we show that it is a
$\BS_{2n}$-submodule of $\bb_n(-2m)$.
\end{abstract} \maketitle

\section{Introduction}

Let $x$ be an indeterminate over $\Z$. The Brauer algebra $\bb_n(x)$
over $\Z[x]$ is a unital associative $\Z[x]$-algebra with generators
$s_1,\cdots,s_{n-1},e_1,$ $\cdots,e_{n-1}$ and relations (see
\cite{E}):
$$\begin{matrix}s_i^2=1,\,\,e_i^2=xe_i,\,\,e_is_i=e_i=s_ie_i,
\quad\forall\,1\leq i\leq n-1,\\
s_is_j=s_js_i,\,\,s_ie_j=e_js_i,\,\,e_ie_j=e_je_i,\quad\forall\,1\leq
i<j-1\leq n-2,\\ s_is_{i+1}s_i=s_{i+1}s_is_{i+1},\,\,
e_ie_{i+1}e_i=e_i,\,\,
e_{i+1}e_ie_{i+1}=e_{i+1},\,\,\forall\,1\leq
i\leq n-2,\\
s_ie_{i+1}e_i=s_{i+1}e_i,\,\,e_{i+1}e_is_{i+1}=e_{i+1}s_i,\quad\forall\,1\leq
i\leq n-2.\end{matrix}
$$
$\bb_n(x)$ is a free $\Z[x]$-module with rank $(2n-1)\cdot
(2n-3)\cdots 3\cdot 1$. For any $\Z[x]$-algebra $R$ with $x$
specialized to $\delta\in R$, we define
$\bb_n(\delta)_{R}:=R\otimes_{\Z[x]}\bb_n(x)$.

This algebra was first introduced by Richard Brauer (see \cite{B})
when he studied how the $n$-tensor space $V^{\otimes{n}}$ decomposes
into irreducible modules over the orthogonal group $O(V)$ or the
symplectic group $Sp(V)$, where $V$ is an orthogonal vector space or
a symplectic vector space. In Brauer's original formulation, the
algebra $\bb_n(x)$ was defined as the complex linear space with
basis the set $\BD_n$ of all the Brauer $n$-diagrams, graphs on $2n$
vertices and $n$ edges with the property that every vertex is
incident to precisely one edge. If we arrange the vertices in two
rows of $n$ each, the top and bottom rows, and label the vertices in
each row of a $n$-diagram by the indices $1,2,\cdots,n$ from left to
right, then the generator $s_i$ corresponds to the $n$-diagram with
edges connecting vertices $i$ (respectively, $i+1$) on the top row
with $i+1$ (respectively, $i$) on bottom row, and all other edges
are vertical, connecting vertex $k$ on the top and bottom rows for
all $k\neq i,i+1$. The generator $e_i$ corresponds to the
$n$-diagram with horizontal edges connecting vertices $i,i+1$ on the
top and bottom rows, and all other edges are vertical, connecting
vertex $k$ on the top and bottom rows for all $k\neq i,i+1$. The
multiplication of two Brauer $n$-diagrams is defined as follows. We
compose two diagrams $D_1, D_2$ by identifying the bottom row of
vertices in the first diagram with the top row of vertices in the
second diagram. The result is a graph, with a certain number,
$n(D_1,D_2)$, of interior loops. After removing the interior loops
and the identified vertices, retaining the edges and remaining
vertices, we obtain a new Brauer $n$-diagram $D_1\circ D_2$, the
composite diagram. Then we define $D_1\cdot
D_2=x^{n(D_1,D_2)}D_1\circ D_2$. In general, the multiplication of
two elements in $\bb_n(x)$ is given by the linear extension of a
product defined on diagrams. For example, let $d$ be the following
Brauer $5$-diagram.
\medskip
\begin{center}
\includegraphics{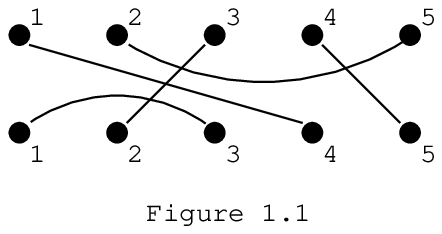}
\end{center}
\medskip
Let $d'$ be the following Brauer $5$-diagram.
\medskip
\begin{center}
\includegraphics{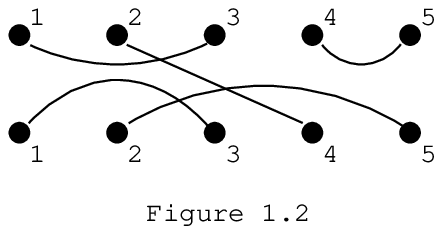}
\end{center}
\medskip
Then $dd'$ is equal to
\medskip
\begin{center}
\includegraphics{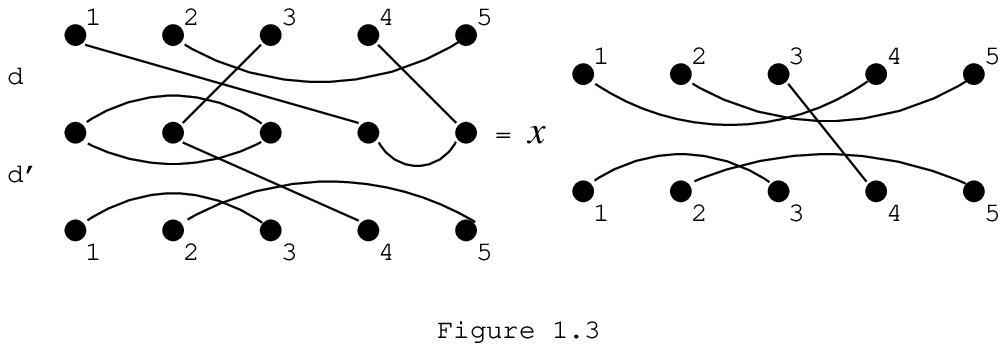}
\end{center}
\medskip
Note that the subalgebra of $\bb_n(x)$ generated by $s_1,s_2,\cdots,s_{n-1}$ is isomorphic to the group algebra
of the symmetric group $\BS_{n}$ over $\Z[x]$.\smallskip

The Brauer algebra as well as its quantization (now called
Birman--Wenzl--Murakami algebra) has been studied in a number of
papers, e.g., \cite{B}, \cite{B1}, \cite{B2}, \cite{HW1},
\cite{HW2}, \cite{BW}, \cite{M}, \cite{Wz}, \cite{FG}, \cite{Xi},
\cite{E}, \cite{DDH}. The walled Brauer algebra (which is a variant
of Brauer algebra, see \cite{BCHLL}) is also studied in the recent
preprint \cite{DD}. We are mainly interested in the Schur--Weyl
duality between symplectic groups and certain specialized Brauer
algebras, which we now recall. Let $K$ be an arbitrary infinite
field. Let $m, n\in\mathbb{N}$. Let $V$ be a $2m$-dimensional
$K$-vector space equipped with a non-degenerate skew-symmetric
bilinear form $(\,,)$. Then (see \cite{Gri}, \cite[Section 4]{Dt})
the symplectic similitude group (respectively, the symplectic group)
relative to $(\,,)$ is
$$ GSp(V):=\Biggl\{g\in GL(V)\Biggm|\begin{matrix}\text{$\exists\, d\in
K$ with $d\neq 0$, such that}\\
\text{$(gv,gw)=d(v,w),\,\,\,\forall\,v,w\in V$}\end{matrix}\Biggr\}
$$
$\Biggl($respectively, $$ Sp(V):=\Bigl\{g\in
GL(V)\Biggm|(gv,gw)=(v,w),\,\,\forall\,\,v,w\in
V\Bigr\}\,\,\Biggr).$$
The symplectic similitude group and
symplectic group $Sp(V)$ act naturally on $V$ from the left-hand
side, and hence on the $n$-tensor space $V^{\otimes n}$. This left
action on $V^{\otimes n}$ is centralized by certain specialized
Brauer algebra, which we recall as follows. Let
$\bb_n(-2m):=\Z\otimes_{\Z[x]}\bb_n(x)$, where $\Z$ is regarded as
$\Z[x]$-algebra by specifying $x$ to $-2m$. Let
$\bb_n(-2m)_K:=K\otimes_{\Z}\bb_n(-2m)$, where $K$ is regarded as
$\Z$-algebra by sending each integer $a$ to $a\cdot 1_{K}$. Then
there is a right action of the specialized Brauer algebra
$\bb_n(-2m)_K$ on the $n$-tensor space $V^{\otimes n}$ which
commutes with the above left action of $GSp(V)$. We recall the
definition of this action as follows. Let $\delta_{ij}$ denote the
value of the usual Kronecker delta. For any $1\leq i\leq 2m$, we set
$$i':=2m+1-i.$$ We fix an ordered basis
$\bigl\{v_1,v_2,\cdots,v_{2m}\bigr\}$ of $V$ such that $$ (v_i,
v_{j})=0=(v_{i'}, v_{j'}),\,\,\,(v_i, v_{j'})=\delta_{ij}=-(v_{j'},
v_{i}),\quad\forall\,\,1\leq i, j\leq m. $$

For any $i, j\in\bigl\{1,2,\cdots,2m\bigr\}$, let $$
\epsilon_{i,j}:=\begin{cases} 1 &\text{if $j=i'$ and $i<j$,}\\
-1 &\text{if $j=i'$ and $i>j$,}\\
0 &\text{otherwise,}\end{cases}
$$
The right action of $\bb_n(-2m)$ on $V^{\otimes n}$ is defined on
generators by $$
\begin{aligned} (v_{i_1}\otimes\cdots\otimes v_{i_n})s_j
&:=-(v_{i_1}\otimes\cdots\otimes v_{i_{j-1}}\otimes
v_{i_{j+1}}\otimes v_{i_{j}}\otimes v_{i_{j+2}}
\otimes\cdots\\
&\qquad\qquad\otimes v_{i_n}),\\
(v_{i_1}\otimes\cdots\otimes v_{i_n})e_j
&:=\epsilon_{i_{j},i_{j+1}} v_{i_1}\otimes\cdots\otimes
v_{i_{j-1}}\otimes \biggl(\sum_{k=1}^{m}(v_{k'}\otimes
v_k-\\
&\qquad\qquad v_{k}\otimes v_{k'})\biggr)\otimes
v_{i_{j+2}}\otimes\cdots\otimes v_{i_n}.\end{aligned}
$$

Let $\varphi,\psi$ be the following natural $K$-algebra
homomorphisms.
$$\begin{aligned}
\varphi:&\,\,(\bb_n(-2m)_K)^{op}\rightarrow\End_{K}\bigl(V^{\otimes
n}\bigr),\\
\psi:&\,\, KGSp(V)\rightarrow\End_{K}\bigl(V^{\otimes n}\bigr)
\end{aligned}
$$

Let $k$ be a positive integer. A composition of $k$ is a sequence of
nonnegative integer $\lam=(\lam_1,\lam_2,\cdots)$ with $\sum_{i\geq
1}\lam_i=k$. A composition $\lam=(\lam_1,\lam_2,\cdots)$ of $k$ is
said to be a partition if $\lam_1\geq\lam_2\geq\cdots$. In this
case, we write $\lam\vdash k$. The conjugate of $\lam$ is defined to
be a partition $\lam'=(\lam'_1,\lam'_2,\cdots)$, where
$\lam'_j:=\#\{i|\lam_i\geq j\}$ for $j=1,2,\cdots$. For any
partition $\lam=(\lam_1,\lam_2,\cdots)$, we use $\ell(\lam)$ to
denote the largest integer $t$ such that $\lam_t\neq 0$.

\begin{lem} {\rm (\cite{B}, \cite{B1}, \cite{B2})} 1) The natural left action of $GSp(V)$ on $V^{\otimes
n}$ commutes with the right action of $\bb_n(-2m)$. Moreover, if
$K=\mathbb{C}$, then
$$\begin{aligned}
\varphi\bigl(\bb_n(-2m)_{\mathbb{C}}\bigr)&=\End_{\mathbb{C}
GSp(V)}\bigl(V^{\otimes n}\bigr)=\End_{\mathbb{C}
Sp(V)}\bigl(V^{\otimes n}\bigr),\\
\psi\bigl(\mathbb{C}GSp(V)\bigr)&=\psi\bigl(\mathbb{C}
Sp(V)\bigr)=\End_{\bb_n(-2m)_{\mathbb{C}}}\bigl(V^{\otimes
n}\bigr),\end{aligned}
$$

2)  if $K=\mathbb{C}$ and $m\geq n$ then $\varphi$ is injective,
and hence an isomorphism onto
$\End_{\mathbb{C}GSp(V)}\bigl(V^{\otimes n}\bigr)$,

3) if $K=\mathbb{C}$, then there is a decomposition of irreducible
$\mathbb{C}GSp(V)$--$\bb_n(-2m)_{\mathbb{C}}$ bimodules
$$
V^{\otimes n}=\bigoplus_{f=0}^{[n/2]}\bigoplus_{\substack{\lam\vdash n-2f\\
\ell(\lam)\leq m}}\Delta({\lam})\otimes D({\lam'}),$$ where
$\Delta({\lam})$ (respectively, $D({\lam'})$) denotes the
irreducible $\mathbb{C}GSp(V)$-module (respectively, the irreducible
$\bb_n(-2m)_{\mathbb{C}}$-module) corresponding to $\lam$
(respectively, corresponding to $\lam'$).
\end{lem}

By recent work of \cite{Oe} and \cite{DDH}, part 1) and part 2) of
the above lemma hold for an arbitrary infinite field. That is,

\begin{lem} {\rm (\cite{Oe}, \cite{DDH})} Let $K$ be an arbitrary infinite field.
\smallskip

1) $\psi\bigl(K GSp(V)\bigr)=\End_{\bb_n(-2m)_K}\bigl(V^{\otimes
n}\bigr).$\smallskip

2) $\varphi\bigl(\bb_n(-2m)_K\bigr)=\End_{KGSp(V)}\bigl(V^{\otimes
n}\bigr) =\End_{KSp(V)}\bigl(V^{\otimes n}\bigr)$, and if $m\geq
n$, then $\varphi$ is also injective, and hence an isomorphism
onto $$
\End_{K GSp(V)}\bigl(V^{\otimes n}\bigr).$$ \end{lem}

Now there is a natural question, that is, how can one describe the
kernel of the homomorphism $\varphi$. This question is closely
related to invariant theory: see \cite{CC}. By the results in
\cite{DDH}, we know that the kernel of the homomorphism $\varphi$
has a rigid structure in the sense that the dimension of
$\Ker\varphi$ does not depend on the choice of the infinite field
$K$, and it is actually defined over $\Z$. Note that in the case of
Schur--Weyl duality between general linear group and symmetric group
(\cite{Sc}, \cite{W}, \cite{CC}, \cite{CL}), or more generally,
between the type $A$ quantum group and the type $A$ Iwahori--Hecke
algebra (\cite{J}, \cite{DPS}), the kernel of the similar
homomorphism has already been explicitly determined in \cite{DPS} in
terms of Kazhdan--Lusztig basis and in \cite{Ha} in terms of Murphy
basis. In this paper, we completely answer the above question by
explicitly constructing an integral basis for the kernel of the
homomorphism $\varphi$. Our description of $\Ker\varphi$ involves a
study of the permutation action of the symmetric group $\BS_{2n}$ on
the Brauer algebra $\bb_n(x)$. Such a permutation action was
previously noted in \cite{FG}. We construct a new integral basis for
this Brauer algebra, which yields integral Specht filtration of this
Brauer algebra by right $\BS_{2n}$-modules. The kernel of $\varphi$
is just one of the $\BS_{2n}$-submodules appearing in this
filtration. In particular, it turns out that $\Ker\varphi$ is in
fact a $\BS_{2n}$-submodule of $\bb_n(-2m)$. The main results of
this paper are presented in Theorem \ref{thm212}, Theorem
\ref{thm214} and Theorem \ref{thm35}. It would be interesting to
compare the new integral basis we obtained in this paper with the
canonical basis for $\bb_n(x)$ constructed in \cite{FG}. Finally, we
remark that it might be possible to give a similar description of
$\Ker\varphi$ also in the orthogonal case (i.e., the case of
Schur--Weyl duality between orthogonal group and certain specialized
Brauer algebra). We deal with only the symplectic case in this paper
because we use the main results in \cite{DDH}, where only the
symplectic case is considered. It would also be interesting to see
how the description of $\Ker\varphi$ we give here can be generalized
to the quantized case, i.e., the case of Schur--Weyl duality between
the quantized enveloping algebra associated to the symplectic Lie
algebra $sp_{2m}$ and a certain specialized Birman--Wenzl--Murakami
algebra (see \cite{CP}).
\bigskip\bigskip

\section{The $\BS_{2n}$-action on $\bb_n(x)$}
\medskip

In this section, we shall first recall (cf. \cite{FG}) the right
permutation action of the symmetric group $\BS_{2n}$ on the set
$\BD_n$. Then we shall construct a new $\Z$-basis for the
resulting right $\BS_{2n}$-module, which yields filtrations of
$\bb_n(x)$ by right $\BS_{2n}$-modules. Certain submodules
occurring in this filtration will play central role in the next
section.
\smallskip

For any involution $\sigma$ in the symmetric group $\BS_{2n}$, the
conjugate $w^{-1}\sigma w$ of $\sigma$ by $w\in\BS_{2n}$ is still
an involution. Therefore, we have a right action of the symmetric
group $\BS_{2n}$ on the set of all the involutions in $\BS_{2n}$.
Note that the set $\BD_n$ of all the Brauer $n$-diagrams can be
naturally identified with the set of all the involutions in
$\BS_{2n}$. Hence we get (cf. \cite{FG}) a right permutation
action of the symmetric group $\BS_{2n}$ on the set $\BD_n$ of all
the Brauer $n$-diagrams. We use $``\ast"$ to denote this right
permutation action. \smallskip

We shall adopt a new labeling of the vertices in each Brauer
diagram. Namely, for each Brauer $n$-diagram $D$, we shall label the
vertices in the top row of $D$ by odd integers $1,3,5,\cdots,2n-1$
from left to right, and label the vertices in the bottom row of $D$
by even integers $2,4,6,\cdots,2n$ from left to right. This way of
labeling is more convenient when studying the permutation action
from $\BS_{2n}$. {\it We shall keep this way of labeling from this
section until the end of Section 3, and we shall recover our
original way of labeling only in Section 4.} Let us look at an
example. Suppose $n=4$, $s_7=(7,8)$ is a transposition in $\BS_{8}$.
Let $D$ be the following Brauer $4$-diagram.

\medskip
\begin{center}
\includegraphics{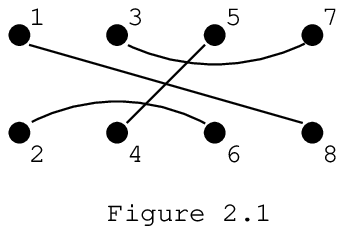}
\end{center}
\medskip

We first identify $D$ with following diagram with $8$-vertices.
\bigskip

\begin{center}
\includegraphics{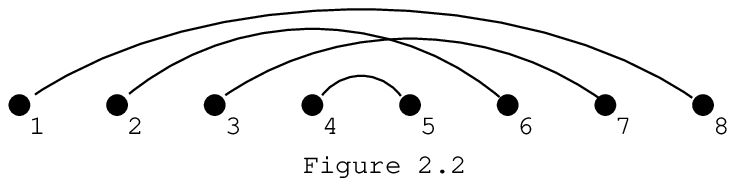}
\end{center}
\medskip

Then $D\ast s_7$ can be computed in the following way.
\bigskip

\begin{center}
\includegraphics{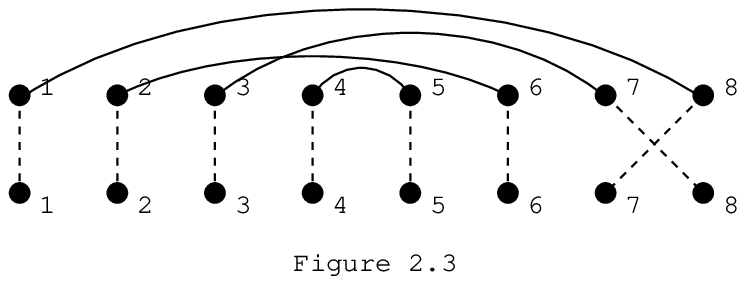}
\end{center}

\medskip
Finally, $D\ast s_7$ is equal to the following Brauer $4$-diagram.
\bigskip

\begin{center}
\includegraphics{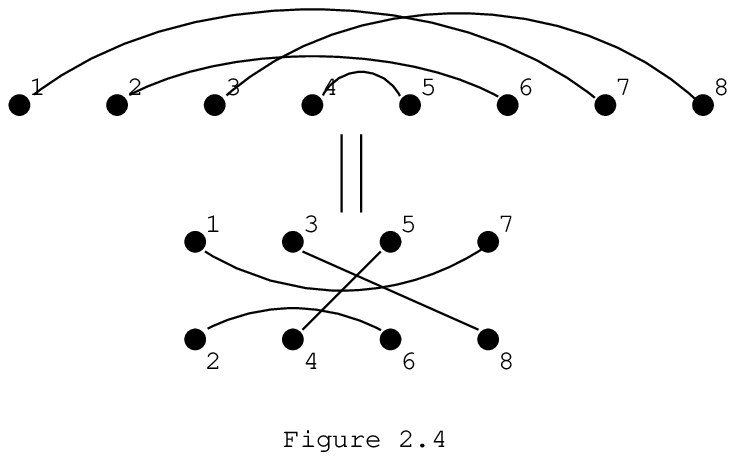}
\end{center}
\medskip

We use $\beta: \BD_n\cong \{w\in\BS_{2n}|w^2=1\}$ to denote the
natural identification of $\BD_n$ with the set of involutions in
$\BS_{2n}$. For any $w\in\BS_{2n}$ and any $D\in\BD_n$, $D\ast
w=\beta^{-1}\bigl(w^{-1}\beta(D)w\bigr)$.
\smallskip

For any $\Z$-algebra $R$, we use $\bb_{n,R}$ to denote the free
$R$-module spanned by all the Brauer $n$-diagrams in $\BD_n$. Then
$\bb_{n,R}$ becomes a right $R[\BS_{2n}]$-module. Let
$\bb_n:=\bb_{n,\Z}$. Clearly, there is a canonical isomorphism
$\bb_{n,R}\cong R\otimes_{\Z}\bb_n$, which is also a right
$R[\BS_{2n}]$-module isomorphism. Taking $R=\Z[x]$, we deduce that
the Brauer algebra $\bb_n(x)$ becomes a right
$\Z[x][\BS_{2n}]$-module. Similarly, the specialized Brauer algebra
$\bb_n(-2m)$ becomes a right $K[\BS_{2n}]$-module.\smallskip

For any integer $i$ with $1\leq i\leq 2n$, we define $$
\gamma({i}):=\begin{cases} i+1, &\text{if $i$ is odd,}\\
i-1, &\text{if $i$ is even.}\end{cases}
$$
Then $\gamma$ is an involution on $\{1,2,\cdots,2n\}$. It is well-known
that the subgroup
$$
\bigl\{w\in\BS_{2n}\bigm|\text{$\bigl(\gamma(a)\bigr)w=\gamma(aw)$
for any integer $a$ with $1\leq a\leq 2n$}\bigr\}
$$
is isomorphic to the wreath product $\Z_2\!\wr\!\BS_n$ of $\Z_2$ and $\BS_n$, which is a Weyl
group of type $B_n$ (c.f. \cite{Hu1}).

\begin{lem} \label{lm21} For any $\Z$-algebra $R$, there is a right $R[\BS_{2n}]$-module isomorphism $$
\bb_{n,R}\cong\Ind_{R[\Z_2\wr\BS_n]}^{R[\BS_{2n}]}1_{R},
$$
where $1_{R}$ denotes the rank one trivial representation of $R[\Z_2\!\wr\!\BS_n]$.
\end{lem}

\begin{proof} Let $1_{\bb_n}$ be the element in $\BD_n$ that connects $2i-1$ to $2i$ for each integer $i$ with
$1\leq i\leq n$. Since $\BS_{2n}$ acts transitively on the set of
all the Brauer $n$-diagrams, it is easy to see that the map
$\xi_{R}$ which send $1_{R}$ to $1_{\bb_n}$ extends naturally to a
surjective $R[\BS_{2n}]$-module homomorphism from
$\Ind_{R[\Z_2\wr\BS_n]}^{R[\BS_{2n}]}1_{R}$ onto $\bb_{n,R}$.

If $R$ is a field, then we can compare the dimensions of both
modules. In that case, we know that the surjection $\xi_R$ must be
an injection, and hence be an isomorphism. In general, since there
are natural isomorphisms
$$ \Ind_{R[\Z_2\wr\BS_n]}^{R[\BS_{2n}]}1_{R}\cong
R\otimes_{\Z}\Ind_{\Z[\Z_2\wr\BS_n]}^{\Z[\BS_{2n}]}1_{\Z},\quad
\bb_{n,R}\cong R\otimes_{\Z}\bb_{n,\Z},
$$
and $\xi_R$ is naturally identified with
$1_{R}\otimes_{\Z}\xi_{\Z}$, it suffices to show that $\xi_{\Z}$ is
an isomorphism. Note also that the short exact sequence
$$ 0\rightarrow \Ker\xi_{\Z}\rightarrow
\Ind_{\Z[\Z_2\wr\BS_n]}^{\Z[\BS_{2n}]}1_{\Z}\rightarrow
\bb_{n,\Z}\rightarrow 0
$$
splits as $\Z$-modules. It follows that $\Ker\xi_{R}$ is canonically
isomorphic to $R\otimes_{\Z}\Ker\xi_{\Z}$ for any $\Z$-algebra $R$.
Let $N:=\Ker\xi_{\Z}$. It is enough to show that $N=0$. By
\cite[Proposition 3.8]{AM}, we only need to show that $N_{(p)}=0$
for each prime number $p$. Let $k_p:=\Z/(p)$, the residue field at
the prime number $p$. It is clear that $k_p\cong
\Z_{(p)}/(p)\Z_{(p)}$. Note that $$ N_{(p)}/(p)N_{(p)}\cong
k_{p}\otimes_{\Z_{(p)}}N_{(p)}\cong
k_{p}\otimes_{\Z_{(p)}}\Ker\xi_{\Z_{(p)}}\cong \Ker\xi_{k_p}=0.
$$
Applying
Nakayama's lemma (\cite[2.6]{AM}), we conclude that $N_{(p)}=0$.
This completes the proof of the lemma.
\end{proof}

For any positive integer $k$ and any composition
$\mu=(\mu_1,\cdots,\mu_s)$ of $k$, the Young diagram of $\mu$ is
defined to be the set $[\mu]:=\{(a,b)|1\leq a\leq s, 1\leq
b\leq\mu_a\}$. The elements of $[\mu]$ are called nodes of $\mu$. A
$\mu$-tableau $\ft$ is defined to be a bijective map from the Young
diagram $[\mu]$ to the set $\{1,2,\cdots,k\}$. For each integer $a$
with $1\leq a\leq k$, we define $\res_{\ft}(a)=j-i$ if $\ft(i,j)=a$.
We denote by $\ft^{\mu}$ the $\mu$-tableau in which the numbers
$1,2,\cdots,k$ appear in order along successive rows. The row
stabilizer of $\ft^{\mu}$, denoted by $\BS_{\mu}$, is the standard
Young subgroup of $\BS_k$ corresponding to $\mu$. For example, if
$k=6, \mu=(2,3,1)$, then $$ \ft^{\mu}=\begin{matrix} 1& 2& \\
3& 4& 5\\
6& &
\end{matrix},\,\,\,\BS_{\mu}=\text{the subgroup of $\BS_6$ generated by $\{s_1,s_3,s_4\}$.}
$$
We define
$$ x_{\mu}=\sum_{w\in\BS_{\mu}}w, \quad
y_{\mu}=\sum_{w\in\BS_{\mu}}(-1)^{\ell(w)}w,
$$
where $\ell(-)$ is the length function in $\BS_k$. If $\mu$ is a
partition of $k$, we denote by $\ft_{\mu}$ the $\mu$-tableau in
which the numbers $1,2,\cdots,k$ appear in order along successive
columns. Let $w_{\mu}\in\BS_{k}$ be such that
$\ft^{\mu}w_{\mu}=\ft_{\mu}$.  For example, if
$k=6, \mu=(3,2,1)$, then $$ \ft^{\mu}=\begin{matrix} 1& 2& 3\\
4& 5\\
6& &
\end{matrix},\,\,\,\,\,
\ft_{\mu}=\begin{matrix} 1& 4& 6\\
2& 5\\
3& &
\end{matrix},\,\,\,\,\,w_{\mu}=(2,4)(3,6).
$$
We use $\mathcal{P}_n$ to denote the set of all the partitions of
$n$. For any partition $\mu$ of $2n$, we define the associated
Specht module $S^{\mu}$ to be the right ideal of the group algebra
$\Z[\BS_{2n}]$ generated by $y_{\mu'}w_{\mu'}x_{\mu}$. In
particular, $S^{(2n)}$ is the one-dimensional trivial representation
of $\BS_{2n}$, while $S^{(1^{2n})}$ is the one dimensional sign
representation of $\BS_{2n}$. By \cite[Theorem 3.5]{DJ2} and
\cite[5.3]{Mu}, our $S^{\mu}$ is isomorphic to the (dual) Specht
module $\widetilde{S}^{\mu}$ introduced in \cite[Section 5]{Mu}. For
any $\Z$-algebra $R$, we write $S_{R}^{\mu}:=R\otimes_{\Z}S^{\mu}$.
Then $\bigl\{S_{\Q}^{\mu}\bigm|\mu\vdash 2n\bigr\}$ is a complete
set of pairwise non-isomorphic simple $\Q[\BS_{2n}]$-modules.
\smallskip

For any composition $\lam=(\lam_1,\cdots,\lam_s)$ of $n$, let
$2\lam:=(2\lam_1,\cdots,2\lam_s)$, which is a composition of $2n$.
We define $2\PP_n:=\bigl\{2\lam\bigm|\lam\in\PP_n\bigr\}$.

\begin{lem} \label{lm22} There is a right $\Q[\BS_{2n}]$-modules isomorphism: $$
\bb_{n,\Q}\cong\bigoplus_{\lam\in 2\PP_n}S_{\Q}^{\lam}.
$$
\end{lem}

\begin{proof} This follows from Lemma \ref{lm21} and \cite[Chapter VII, (2.4)]{Mac}.
\end{proof}

Let $a$ be an integer with $0\leq a\leq n$. Let $1\leq
i_1,\cdots,i_a, j_1,\cdots,j_a\leq 2n$ be $2a$ pairwise distinct
integers. Let
$$I:=\{1,2,\cdots,2n\}\setminus\{i_1,\cdots,i_a,j_1,\cdots,j_a\}.$$
Let $\BS_{I}$ be the symmetric group on the set $I$. Let
$\ui:=(i_1,\cdots,i_a), \uj:=(j_1,\cdots,j_a)$. Let $$
\BD_n(\ui,\uj):=\Bigl\{D\in\BD_n\Bigm|\text{$D$ connects $i_s$
with $j_s$ for each $1\leq s\leq a$}\Bigr\}.
$$

\begin{lem} \label{lm23} With the notations as above, for any $w\in\BS_{I}$, we have $$
\Biggl(\sum_{D\in\BD_n(\ui,\uj)}D\Biggr)\ast
w=\sum_{D\in\BD_n(\ui,\uj)}D.
$$
\end{lem}

\begin{proof} For any $D\neq D'\in\BD_n(\ui,\uj)$, it is clear that
$$D\ast w\neq D'\ast w\in\BD_n(\ui,\uj).$$ Therefore, the lemma follows easily from a counting argument.
\end{proof}

\addtocounter{dfn}{3}
\begin{dfn} For any non-negative even integers $a,b$ with $a+b\leq 2n$, we define $$
\begin{aligned}
\BD_{(b)}^{(a)}:&=\Biggl\{D\in\BD_n\Biggm|\begin{matrix}\text{the vertex labeled by $i$ is connected with}\\
\text{the vertex labeled by $\theta(i)$ whenever}\\
\text{$i\leq a$ or $i>a+b$}\end{matrix}\Biggr\},\\
X_{(b)}^{(a)}:&=\sum_{D\in\BD_{(b)}^{(a)}}D.
\end{aligned}
$$
\end{dfn}

For any even integer $k$ with $0\leq k\leq 2n$, let $X_{(k)}:=X_{(k)}^{(0)}$.\smallskip

\begin{dfn} Let $\lam=(\lam_1,\cdots,\lam_s)$ be a composition of $2n$ such that $\lam_i$ is even for each $i$,
we define $$
X_{\lam}:=X_{(\lam_1)}^{(0)}X_{(\lam_2)}^{(\lam_1)}\cdots X_{(\lam_s)}^{(\lam_1+\cdots+\lam_{s-1})}\in\bb_n.
$$
\end{dfn}

\addtocounter{cor}{5}
\begin{cor} \label{cor26} Let $\lam=(\lam_1,\cdots,\lam_s)$ be a composition of $2n$ such that $\lam_i$ is even
for each $i$. Then, for any $w\in\BS_{\lam}$, we have that
$$ X_{\lam}\ast w=X_{\lam}.
$$
\end{cor}

\begin{proof} For any non-negative even integers $a,b$ with $a+b\leq 2n$, the set $\BD_{(b)}^{(a)}$ is just a special
case of the set $\BD_n(\ui,\uj)$ we defined before. Therefore, by
Lemma \ref{lm23}, for any $w\in\BS_{\{a+1,a+2,\cdots,a+b\}}$, we
have $$ X_{(b)}^{(a)}\ast w=X_{(b)}^{(a)}.
$$
Now we note that the elements $X_{(\lam_1)}^{(0)}, X_{(\lam_2)}^{(\lam_1)}, \cdots,
X_{(\lam_s)}^{(\lam_1+\cdots+\lam_{s-1})}$ pairwise commute with each other. Hence the corollary follows at once.
\end{proof}

Let $k$ be a positive integer and $\mu$ be a composition of $k$. A
$\mu$-tableau $\ft$ is called {\it row standard} if the numbers
increase along rows. We use $\RS(\mu)$ to denote the set of all the
row-standard $\mu$-tableaux. Suppose $\mu$ is a partition of $k$.
Then $\ft$ is called {\it column standard} if the numbers increase
down columns, and {\it standard} if it is both row and column
standard. In this case, it is clear that both $\ft^{\mu}$ and
$\ft_{\mu}$ are standard $\mu$-tableaux. We use $\Std(\mu)$ to
denote the set of all the standard $\mu$-tableaux. Now let $\mu$ be
a composition of $2n$. Let $\lam\in 2\PP_n$. For any
$\ft\in\RS(\lam)$, let $d(\ft)\in\BS_{2n}$ be such that
$\ft^{\lam}d(\ft)=\ft$. Let $X_{\lam,\ft}:=X_{\lam}\ast d(\ft)$. For
any $\Z$-algebra $R$, we define $$
\CM^{\lam}_{R}:=\text{$R$-$\Span$}\Bigl\{X_{\nu,\ft}\Bigm|\ft\in\Std(\nu),\lam\unlhd\nu\in
2\PP_n\Bigr\},
$$
where ``$\unrhd$'' is the dominance order defined in \cite{Mu}. We
write $\CM^{\lam}=\CM_{\Z}^{\lam}$. We are interested in the module
$\CM_{R}^{\lam}$. In the remaining part of this paper, we shall see
that this module is actually a right $\BS_{2n}$-submodule of
$\bb_{n,R}$, and it shares many properties with the permutation
module $x_{\lam}\Z[\BS_{2n}]$. In particular, it also has a Specht
filtration, and it is stable under base change, i.e.,
$R\otimes_{\Z}\CM^{\lam}\cong \CM_{R}^{\lam}$ for any $\Z$-algebra
$R$.\smallskip

For our purpose, we need to recall some results in \cite{Mu} and
\cite{Ma} on the Specht filtrations of permutation modules over the
symmetric group $\BS_{2n}$. Let $\lam$, $\mu$ be two partitions of
$2n$. A $\mu$-tableau of type $\lam$ is a map
$\fS:[\mu]\rightarrow\{1,2,\cdots,2n\}$ such that each $i$ appears
exactly $\lam_i$ times. $\fS$ is said to be semistandard if each row
of $\fS$ is weakly increasing and each column of $\fS$ is strictly
increasing. Let $\T_0(\mu,\lam)$ be the set of all the semistandard
$\mu$-tableaux of type $\lam$. Then $\T_0(\mu,\lam)\neq\emptyset$
only if $\mu\unrhd \lam$. For each standard $\mu$-tableau $\fs$, let
$\mu(\fs)$ be the tableau which is obtained from $\fs$ by replacing
each entry $i$ in $\fs$ by $r$ if $i$ appears in row $r$ of
$\ft^{\lam}$. Then $\mu(\fs)$ is a $\mu$-tableau of type
$\lam$.\smallskip

For each standard $\mu$-tableau $\ft$ and each semistandard $\mu$-tableau $\fS$ of type $\lam$, we define $$
x_{\fS,\ft}:=\sum_{\fs\in\Std(\mu),\mu(\fs)=\fS}d(\fs)^{-1}x_{\mu}d(\ft).
$$
Then by \cite[Section 7]{Mu}, the set $$
\Bigl\{x_{\fS,\ft}\Bigm|\fS\in\T_0(\mu,\lam), \ft\in\Std(\mu), \lam\unlhd\mu\vdash 2n\Bigr\}
$$
form a $\Z$-basis of $x_{\lam}\Z[\BS_{2n}]$. Furthermore, for any
$\Z$-algebra $R$, the canonical surjective homomorphism
$R\otimes_{\Z}x_{\lam}\Z[\BS_{2n}]\twoheadrightarrow
x_{\lam}R[\BS_{2n}]$ is an isomorphism.

For each partition $\mu$ of $2n$ and for each semistandard
$\mu$-tableau $\fS$ of type $\lam$, according to the results in
\cite[Section 7]{Mu} and \cite{Ma}, both the following
$\Z$-submodules $$\begin{aligned}
M^{\lam}_{\fS}:&=\text{$\Z$-$\Span$}\Bigl\{x_{\fS,\fs},
x_{\fT,\ft}\Bigm|\begin{matrix}\fs\in\Std(\mu),
\fT\in\T_0(\nu,\lam),\\ \ft\in\Std(\nu),\mu\lhd\nu\vdash 2n\end{matrix}\Bigr\},\\
M^{\lam}_{\fS,\rhd}:&=\text{$\Z$-$\Span$}\Bigl\{x_{\fT,\ft}\Bigm|\fT\in\T_0(\nu,\lam),
\ft\in\Std(\nu),\mu\lhd\nu\vdash 2n\Bigr\},
\end{aligned}
$$
are $\Z[\BS_{2n}]$-submodules, and the quotient of $M^{\lam}_{\fS}$
by $M^{\lam}_{\fS,\rhd}$ is canonically isomorphic to $S^{\mu}$ so
that the image of the elements $x_{\fS,\fs}$, where
$\fs\in\Std(\mu)$, forms the standard $\Z$-basis of $S^{\mu}$. In
other words, it gives rise to the Specht filtrations of
$x_{\lam}\Z[\BS_{2n}]$, each semistandard $\mu$-tableau of type
$\lam$ yields a factor which is isomorphic to $S^{\mu}$ so that
$x_{\lam}\Z[\BS_{2n}]$ has a series of factors, ordered by $\unlhd$,
each isomorphic to some $S^{\mu}$, $\mu\unrhd \lam$, the
multiplicity of $S^{\mu}$ being the number of semistandard
$\mu$-tableaux of type $\lam$.

We write $\mu=(\mu_1,\mu_2,\cdots)=(a_1^{k_1},a_2^{k_2},\cdots,)$,
where $a_1>a_2>\cdots$, $k_i\in\mathbb{N}$ for each $i$, where
$a_i^{k_i}$ means that $a_i$ repeats $k_i$ times. Let
$\widetilde{\BS}_{\mu}$ be the subgroup of $\BS_{\mu'}$ consisting
of all the elements $w$ satisfying the following condition: for any
integer $t\geq 1$, and any integers $i, j$ with
$\sum_{s=1}^{t-1}k_s+1\leq i, j\leq\sum_{s=1}^{t}k_s$, and any
integers $a,b$ with $1\leq a, b\leq\mu_{\sum_{s=1}^{t}k_s}$,
\addtocounter{equation}{6} \begin{equation}
\text{$(\ft_{\mu}(i,a))w=\ft_{\mu}(j,a)$\quad if and only if \quad
$(\ft_{\mu}(i,b))w=\ft_{\mu}(j,b)$}. \end{equation} Let
$\widetilde{D}_{\mu}$ be a complete set of right coset
representatives of $\widetilde{\BS}_{\mu}$ in $\BS_{\mu'}$.

\addtocounter{lem}{4}
\begin{lem} \label{lm27} For any partition $\lam\in 2\PP_n$, let $$
n_{\lam}:=\prod_{i\geq 1}(\lam_i-\lam_{i+1})!,\quad
h_{\lam}:=\sum_{w\in\widetilde{D}_{\lam}}(-1)^{\ell(w)}w.
$$
Then $$
X_{\lam}\ast \bigl(w_{\lam}y_{\lam'}\bigr)=n_{\lam}\bigl(X_{\lam}\ast (w_{\lam}h_{\lam})\bigr),
$$
and for any $\Z$-algebra $R$, $1_{R}\otimes_{\Z}(X_{\lam}\ast (w_{\lam}h_{\lam}))\neq 0$ in
$\bb_{n,R}$.
\end{lem}

\begin{proof} The condition $\lam\in 2\PP_n$ implies that for any $w\in\widetilde{\BS}_{\lam}$, $\ell(w)$ is an even integer. The first statement of this lemma now follows from the following identity: $$
(X_{\lam}\ast w_{\lam})\ast \Bigl(\sum_{w\in\widetilde{\BS}_{\lam}}w\Bigr)=n_{\lam}X_{\lam}\ast w_{\lam}.
$$
Let $d$ be the Brauer $n$-diagram in which the vertex labeled by
$\ft_{\lam}(i,2j-1)$ is connected with the vertex labeled by
$\ft_{\lam}(i,2j)$ for any $1\leq i\leq\lam'_1, 1\leq
j\leq\lam_i/2$. Then it is easy to see that $d$ appears with
coefficient $1$ in the expression of $X_{\lam}\ast
(w_{\lam}h_{\lam})$ as linear combinations of basis of Brauer
$n$-diagrams. It follows that for any $\Z$-algebra $R$,
$1_{R}\otimes_{\Z}(X_{\lam}\ast (w_{\lam}h_{\lam}))\neq 0$ in
$\bb_{n,R}$, as required.
\end{proof}

Following \cite{Mu0}, we define the Jucys-Murphy operators of
$\Z[\BS_{2n}]$.
$$\left\{\begin{aligned}
L_1:&=0,\\
L_a:&=(a-1,a)+(a-2,a)+\cdots+(1,a),\quad a=2,3,\cdots,2n.
\end{aligned}
\right.$$ Then for each partition $\lam$ of $2n$, and each integer
$1\leq a\leq 2n$, we have (by \cite[(3.14)]{DJ2}) $$
\bigl(x_{\lam}w_{\lam}y_{\lam'}\bigr)L_a=\res_{\ft_{\lam}}(a)\bigl(x_{\lam}w_{\lam}y_{\lam'}\bigr).
$$
For each standard $\lam$-tableau $\ft$, we define $$
\Theta_{\ft}:=\prod_{i=1}^{n}\prod_{\substack{\fu\in\Std(\lam)\\
\res_{\fu}(i)\neq\res_{\ft}(i)}}\frac{L_i-\res_{\fu}(i)}{\res_{\ft}(i)-\res_{\fu}(i)}.
$$
For each partition $\lam\in 2\PP_n$, by Corollary \ref{cor26} and
Frobenius reciprocity, there is a surjective  right
$\Z[\BS_{2n}]$-module homomorphism $\pi_{\lam}$ from
$x_{\lam}\Z[\BS_{2n}]$ onto $X_{\lam}\Z[\BS_{2n}]$ which extends
the map $x_{\lam}\mapsto X_{\lam}$. In particular, by Lemma \ref{lm27}, $$
\bigl(X_{\lam}\ast w_{\lam}h_{\lam}\bigr)\ast L_a=\res_{\ft_{\lam}}(a)\bigl(X_{\lam}\ast (w_{\lam}h_{\lam})\bigr).
$$

\addtocounter{prop}{8}
\begin{prop} \label{prop29} For any partition $\lam\in 2\PP_n$, we have that
$$[X_{\lam}\Q[\BS_{2n}]:S_{\Q}^{\lam}]=1. $$
\end{prop}

\begin{proof} By Lemma \ref{lm22}, we have $$
\bb_{n,\Q}\cong\bigoplus_{\mu\in 2\PP_n}S_{\Q}^{\mu}.
$$
It is well-known that each $S_{\Q}^{\mu}$ has a basis
$\bigl\{v_{\ft}\bigr\}_{\ft\in\Std(\mu)}$ satisfying $$
v_{\ft}L_i=\res_{\ft}(i)v_{\ft},\quad \forall\,1\leq i\leq n.
$$
Let $\lam$ be a fixed partition in $2\PP_n$. Since $X_{\lam}\Q[\BS_{2n}]\subseteq \bb_{n,\Q}$, we can write $$
X_{\lam}\ast (w_{\lam}h_{\lam})=\sum_{\mu\in
2\PP_n}\sum_{\ft\in\Std(\mu)}A_{\ft}v_{\ft},
$$
where $A_{\ft}\in\Q$ for each $\ft$.

For each $\mu\in 2\PP_n$ and each $\ft\in\Std(\mu)$, we apply
the operator $\Theta_{\ft}$ on both sides of the above identity
and use Lemma \ref{lm27} and the above discussion. We get that
$A_{\ft}\neq 0$ if and only if $\mu=\lam$ and $\ft=\ft_{\lam}$. In other words,
$X_{\lam}\ast (w_{\lam}h_{\lam})=A_{\ft_{\lam}}v_{\ft_{\lam}}$ for
some $0\neq A_{\ft_{\lam}}\in\Q$. This implies that the projection
from $X_{\lam}\Q[\BS_{2n}]$ to $S_{\Q}^{\lam}$ is nonzero.
Hence, $[X_{\lam}\Q[\BS_{2n}]:S_{\Q}^{\lam}]=1$, as required.
\end{proof}

For each partition $\lam\in 2\PP_n$, by the natural surjective
$\Z[\BS_{2n}]$-module homomorphism $\pi_{\lam}$ from
$x_{\lam}\Z[\BS_{2n}]$ onto $X_{\lam}\Z[\BS_{2n}]$, we know that the
elements $\pi_{\lam}\bigl(x_{\fS,\ft}\bigr)$, where
$\fS\in\T_0(\mu,\lam), \ft\in\Std(\mu), \lam\unlhd\mu\vdash 2n$,
spans $X_{\lam}\Z[\BS_{2n}]$ as $\Z$-module.

\begin{prop} \label{prop210} For any two partitions $\lam, \mu$ of $2n$, and for any $\fS\in\T_0(\mu,\lam)$, we have
that $\pi_{\lam}\bigl(M_{\fS}^{\lam}\bigr)\subseteq\CM^{\lam}$. In particular,
$X_{\lam}\Z[\BS_{2n}]\subseteq \CM^{\lam}$.
\end{prop}

\begin{proof} We first prove a weak version of the claim in this
proposition. That is, for any two partitions $\lam, \mu$ of $2n$,
and for any $\fS\in\T_0(\mu,\lam)$, $$
\pi_{\lam}\bigl(M_{\fS}^{\lam}\bigr)\subseteq\CM_{\Q}^{\lam}.$$

We consider the dominance order $``\unlhd"$ and make induction on
$\lam$. We start with the partition $(2n)$, which is the unique
maximal partition of $2n$ with respect to $``\unlhd"$. In this
case, $x_{(2n)}\Z[\BS_{2n}]=\Z x_{(2n)}$, and
$X_{(2n)}\Z[\BS_{2n}]=\Z X_{(2n)}$, it is easy to see the claim in
this proposition is true for $\lam=(2n)$.

Now let $\lam\lhd (2n)$ be a partition of $2n$. Assume that for any
partition $\nu$ of $2n$ satisfying $\nu\rhd\lam$, the claim in
this proposition is true. We now prove the claim for the partition
$\lam$.

Let $\mu\unrhd\lam$ be a partition of $2n$ with
$\T_0(\mu,\lam)\neq\emptyset$. We consider again the dominance order
$``\unlhd"$ and make induction on $\mu$. Since $\T_0((2n),\lam)$
contains a unique element $\fS_{\star}$, $\Std((2n))=\{\ft^{(2n)}\}$
and
$$
\pi_{\lam}\bigl(x_{\fS_{\star},\ft^{(2n)}}\bigr)=\pi_{\lam}(x_{(2n)})
=X_{(2n)}\in\CM^{\lam}.$$
So in this case the claim of this proposition is still true.

Now let $\mu\unrhd\lam$ be a partition of $2n$ with
$\T_0(\mu,\lam)\neq\emptyset$ and $\mu\lhd (2n)$. Assume that for
any partition $\nu$ of $2n$ satisfying $\T_0(\nu,\lam)\neq\emptyset$
and $\nu\rhd\mu$, and for any $S\in\T_0(\nu,\lam)$, $$
\pi_{\lam}\bigl(M_{\fS}^{\lam}\bigr)\subseteq\CM_{\Q}^{\lam}.
$$
Let $\fS\in\T_0(\mu,\lam)$. The homomorphism $\pi_{\lam}$ induces a
surjective map from $M_{\fS}^{\lam}/M_{\fS,\rhd}^{\lam}$ onto $$
\Bigl(\pi_{\lam}(M_{\fS}^{\lam})\Bigr)/\Bigl(\pi_{\lam}(M_{\fS,\rhd}^{\lam})\Bigr).
$$
Hence it also induces a surjective map $\widetilde{\pi}_{\lam}$ from $$
\Bigl(\Q\otimes_{\Z}M_{\fS}^{\lam}/\Bigl(\Q\otimes_{\Z}M_{\fS,\rhd}^{\lam}\Bigr)\cong
\Q\otimes_{\Z}\Bigl(M_{\fS}^{\lam}/M_{\fS,\rhd}^{\lam}\Bigr)\cong S_{\Q}^{\mu}
$$
onto $$
\Q\otimes_{\Z}\Bigl(\pi_{\lam}(M_{\fS}^{\lam})/\pi_{\lam}(M_{\fS,\rhd}^{\lam})\Bigr).
$$
Since $S_{\Q}^{\mu}$ is irreducible, the above map is either a zero
map or an isomorphism. If it is a zero map, then (by induction
hypothesis) $$
\pi_{\lam}(M_{\fS}^{\lam})\subseteq\pi_{\lam}(M_{\fS,\rhd}^{\lam})_{\Q}\subseteq
\CM_{\Q}^{\lam}.
$$
It remains to consider the case where $\widetilde{\pi}_{\lam}$ is an isomorphism.
In particular, $$\Q\otimes_{\Z}\Bigl(\pi_{\lam}(M_{\fS}^{\lam})/\pi_{\lam}(M_{\fS,\rhd}^{\lam})\Bigr)\cong S_{\Q}^{\mu}.
$$
Applying Lemma \ref{lm22}, we know that $\mu\in 2\PP_n$.\smallskip

On the other hand, the homomorphism $\pi_{\mu}$ also induces a surjective map from
$x_{\mu}\Z[\BS_{2n}]/M_{\fS_0,\rhd}^{\mu}$
onto $$
\Bigl(\pi_{\mu}(x_{\mu}\Z[\BS_{2n}])\Bigr)/\Bigl(\pi_{\mu}(M_{\fS_0,\rhd}^{\mu})\Bigr)=
X_{\mu}\Z[\BS_{2n}]/\Bigl(\pi_{\mu}(M_{\fS_0,\rhd}^{\mu})\Bigr),
$$
where $\fS_0$ is the unique semistandard $\mu$-tableau in $\T_0(\mu,\mu)$.
Hence it also induces a surjective map $\widetilde{\pi}_{\mu}$ from $$
\Bigl(\Q\otimes_{\Z}x_{\mu}\Z[\BS_{2n}]/\Bigl(\Q\otimes_{\Z}M_{\fS_0,\rhd}^{\mu}\Bigr)\cong
\Q\otimes_{\Z}\Bigl(x_{\mu}\Z[\BS_{2n}]/M_{\fS_0,\rhd}^{\mu}\Bigr)\cong S_{\Q}^{\mu}
$$
onto $$
\Q\otimes_{\Z}\Bigl(X_{\mu}\Z[\BS_{2n}]/\pi_{\mu}(M_{\fS_0,\rhd}^{\mu})\Bigr)\cong
\Bigl(X_{\mu}\Q[\BS_{2n}]\Bigr)/\Bigl(\Q\otimes_{\Z}\pi_{\mu}(M_{\fS_0,\rhd}^{\mu})\Bigr).
$$

By the Specht filtration of $M_{\Q}^{\mu}$, we know that
$S_{\Q}^{\mu}$ does not occur as composition factor in
$\Q\otimes_{\Z}M_{\fS_0,\rhd}^{\mu}$. Hence $S_{\Q}^{\mu}$ does not
occur as composition factor in
$\Q\otimes_{\Z}\pi_{\mu}(M_{\fS_0,\rhd}^{\mu})$. By Proposition
\ref{prop29}, $S_{\Q}^{\mu}$ occurs as composition factor with
multiplicity one in $X_{\mu}\Q[\BS_{2n}]$. Therefore,
$X_{\mu}\Q[\BS_{2n}]\neq\Q\otimes_{\Z}\pi_{\mu}(M_{\fS_0,\rhd}^{\mu})$.
It follows that $\widetilde{\pi}_{\mu}$ must be an isomorphism.
Hence
$$
\Q\otimes_{\Z}\Bigl(X_{\mu}\Z[\BS_{2n}]/\pi_{\mu}(M_{\fS_0,\rhd}^{\mu})\Bigr)\cong
S_{\Q}^{\mu}.
$$

We write $A=\pi_{\lam}(M_{\fS}^{\lam}), B=X_{\mu}\Z[\BS_{2n}]$. Since $S_{\Q}^{\mu}$ appears only once
in $\bb_{n,\Q}$, it follows that $S_{\Q}^{\mu}$ must occur as composition factor in the module $$
\bigl(\Q\otimes_{\Z}A\bigr)\cap \bigl(\Q\otimes_{\Z}B\bigr)=\Q\otimes_{\Z}(A\cap B).
$$
Hence $S_{\Q}^{\mu}$ can not occur as composition factor in the module $$
\bigl(\Q\otimes_{\Z}A\bigr)/\Bigl(\Q\otimes_{\Z}(A\cap B)\Bigr)\cong\Q\otimes_{\Z}(A/A\cap B).
$$
Therefore, the image of the canonical projection
$\Q\otimes_{\Z}A\rightarrow\Q\otimes_{\Z}(A/A\cap B)$ must  be
contained in the image of
$\Q\otimes_{\Z}\pi_{\lam}(M_{\fS,\rhd}^{\lam})$. However, by
induction hypothesis, both $\pi_{\lam}(M_{\fS,\rhd}^{\lam})$ and
$B$ are contained in the $\Q$-span of
$\Bigl\{X_{\alpha,\fu}\Bigm|\fu\in\Std(\alpha),
\lam\unlhd\alpha\in 2\PP_n\Bigr\}$. It follows that $$
\pi_{\lam}\bigl(M_{\fS}^{\lam}\bigr)\subseteq\CM_{\Q}^{\lam},
$$
as required.\smallskip

Now we begin to prove
$\pi_{\lam}\bigl(M_{\fS}^{\lam}\bigr)\subseteq\CM^{\lam}$. Suppose
that $$ \pi_{\lam}\bigl(M_{\fS}^{\lam}\bigr)\not\subseteq\CM^{\lam}.
$$
Then (by the $\Z$-freeness of $\bb_n$) there exist an element $x\in
M_{\fS}^{\lam}$, some integers $a, a_{\fu}$, and a prime divisor
$p\in\mathbb{N}$ of $a$, such that $$
a\pi_{\lam}(x)=\sum_{\lam\unlhd\alpha\in
2\PP_n}\sum_{\fu\in\Std(\alpha)}a_{\fu}X_{\alpha}\ast d(\fu),
$$
and $$\Sigma_p:=\bigl\{\alpha\in
2\PP_n\bigm|\text{$\lam\unlhd\alpha, p\nmid a_{\fu}$, for some
$\fu\in\Std(\alpha)$}\bigr\}\neq\emptyset. $$

We take an $\alpha\in\Sigma_p$ such that $\alpha$ is minimal with respect to ``$\unlhd$". Then we take an
$u\in\Std(\alpha)$ such that $p\nmid a_{\fu}$ and $\ell(d(\fu))$ is maximal among the elements in the set
$\bigl\{\fu\in\Std(\alpha)\bigm|p\nmid a_{\fu}\bigr\}$. Let $\sigma_{\fu}$ be the unique element in $\BS_{2n}$ such that
$d(\fu)\sigma_{\fu}=w_{\alpha}$ and $\ell(w_{\alpha})=\ell(d(\fu))+\ell(\sigma_{\fu})$. We consider the finite field $\mathbb{F}_p$ as a $\Z$-algebra. By \cite[(4.1)]{DJ1}, we know that for any partitions $\beta, \gamma$ of $2n$, and element $w\in\BS_{2n}$, $$
\text{$x_{\beta}wy_{\gamma'}\neq 0$ only if $\gamma\unrhd\beta$};\,\,\text{while}\,\,\text{$x_{\beta}wy_{\beta'}\neq 0$ only if $w\in\BS_{\beta}w_{\beta}$}.
$$
Hence by Lemma \ref{lm27} and the homomorphism $\pi_{\lam}$, $$
\text{$X_{\beta}\ast (wh_{\gamma})\neq 0$ only if
$\gamma\unrhd\beta$};\,\,\text{$X_{\beta}\ast (wh_{\beta})\neq 0$
only if $w\in\BS_{\beta}w_{\beta}$}. $$ Using Lemma \ref{lm27}
again, we get $$
0=1_{\mathbb{F}_p}\otimes_{\Z}\bigl(a\pi_{\lam}(x)\ast
(\sigma_{\fu}h_{\alpha})\bigr)=1_{\mathbb{F}_p}\otimes_{\Z}\bigl(a_{\fu}X_{\alpha}\ast
(w_{\alpha}h_{\alpha})\bigr)\neq 0,
$$
which is a contradiction. This prove that $\pi_{\lam}\bigl(M_{\fS}^{\lam}\bigr)\subseteq\CM^{\lam}$.
\end{proof}

\addtocounter{cor}{4}
\begin{cor} \label{cor211} For any partition $\lam\in 2\PP_n$ and any $\Z$-algebra $R$, $\CM_{R}^{\lam}$ is a right $\BS_{2n}$-submodule of $\bb_{n,R}$. \end{cor}

\begin{proof} This follows directly from Proposition \ref{prop210}.
\end{proof}

\addtocounter{thm}{11}
\begin{thm} \label{thm212} For any partition $\lam\in 2\PP_n$ and any $\Z$-algebra $R$, the canonical map
$R\otimes_{\Z}\CM^{\lam}\rightarrow\CM_{R}^{\lam}$ is an isomorphism, and the set $$
\Bigl\{X_{\nu,\ft}\Bigm|\ft\in\Std(\nu),\lam\unlhd\nu\in 2\PP_n\Bigr\}
$$
form an $R$-basis of $\CM_{R}^{\lam}$. In particular, the set $$
\Bigl\{X_{\lam,\ft}\Bigm|\ft\in\Std(\lam),\lam\in 2\PP_n\Bigr\}
$$
form an $R$-basis of $\bb_{n,R}$.
\end{thm}

\begin{proof} We take $\lam=(2^{n})$, then $X_{\lam}\Z[\BS_{2n}]=\bb_n$. It is well-known that $
\bb_{n,R}\cong R\otimes_{\Z}\bb_{n}$ for any $\Z$-algebra $R$.
Applying Proposition \ref{prop210} and counting the dimension, we
get that for any $\Z$-algebra $R$ which is a field, the set $$
\Bigl\{X_{\lam,\ft}\Bigm|\ft\in\Std(\lam),\lam\in 2\PP_n\Bigr\}
$$
must form an $R$-basis of $\bb_{n,R}$. Since each basis element
$X_{\lam,\ft}$ is integrally defined and $\bb_{n,R}\cong
R\otimes_{\Z}\bb_{n}$ for any $\Z$-algebra $R$, it follows that the
above set is still an $R$-basis of $\bb_{n,R}$ for any $\Z$-algebra
$R$.\smallskip

By the $R$-linear independence of the elements in this set and
Corollary \ref{cor211}, we also get that, for any partition $\lam\in
2\PP_n$, the set $$
\Bigl\{X_{\nu,\ft}\Bigm|\ft\in\Std(\nu),\lam\unlhd\nu\in
2\PP_n\Bigr\}
$$
must form an $R$-basis of $\CM_{R}^{\lam}$. Therefore, for any $\Z$-algebra $R$, the canonical map $R\otimes_{\Z}\CM^{\lam}\rightarrow \CM_{R}^{\lam}$ is an isomorphism.
\end{proof}

\noindent
{\it Remark 2.13.} Note that for any partition $\lam\in 2\PP_n$, $X_{\lam}\Z[\BS_{2n}]\subseteq\CM^{\lam}$. But $X_{\lam}\Z[\BS_{2n}]$ is not necessarily equal to $\CM^{\lam}$ in general. For example, one sees easily that $$
X_{(6,2)}\notin X_{(4,4)}\Z[\BS_{8}].
$$
In fact, if this is not the case, then we can write $$
X_{(6,2)}+\sum_{i}a_iX_{(4,4)}\ast w_i=\sum_{i}b_jX_{(4,4)}\ast w'_j,
$$
for some positive integers $a_i, b_j$ and some elements $w_i,
w'_j\in\BS_{2n}$. However, if we express both sides into linear
combinations of Brauer $4$-diagrams and count the number of terms,
we find that this is impossible (as the equation $15+9a=9b$ has no
solutions in $\Z$). \smallskip

\addtocounter{thm}{1}
\begin{thm} \label{thm214} For any partition $\lam\in 2\PP_n$ and any $\Z$-algebra $R$, we define $$
\CM_R^{\rhd\lam}:=\text{$R$-$\Span$}\Bigl\{X_{\nu,\ft}\Bigm|\ft\in\Std(\nu),\lam\lhd\nu\in 2\PP_n\Bigr\}.
$$
Then $\CM_R^{\rhd\lam}$ is a right $R[\BS_{2n}]$-submodule of $\CM_R^{\lam}$, and there is a $R[\BS_{2n}]$-module
isomorphism $$
\CM_R^{\lam}/\CM_R^{\rhd\lam}\cong S_R^{\lam}.
$$
In particular, $\bb_{n,R}$ has a Specht filtration.
\end{thm}

\begin{proof} It suffices to consider the case where $R=\Z$. We first show that $$
\CM_{\Q}^{\lam}\cong\oplus_{\lam\unlhd\mu\in 2\PP_n}S_{\Q}^{\mu},\quad
\CM_{\Q}^{\rhd\lam}\cong\oplus_{\lam\lhd\mu\in 2\PP_n}S_{\Q}^{\mu}.
$$
For each $\mu\in 2\PP_n$, we use $\rho_{\mu}^{\lam}$ to denote the composition of the embedding $\CM_{\Q}^{\lam}\hookrightarrow
\bb_{n,\Q}$ and the projection $\bb_{n,\Q}\twoheadrightarrow S_{\Q}^{\mu}$. Suppose that $\rho_{\mu}^{\lam}\neq 0$. Then
$\rho_{\mu}^{\lam}$ must be an surjection. We claim that $\mu\unrhd\lam$. In fact, if $\mu\ntrianglerighteq\lam$, then
for any $\lam\unlhd\nu\in 2\PP_n$, $\mu\ntrianglerighteq\nu$, and $x_{\nu}\Z[\BS_{2n}]w_{\mu'}x_{\mu}w_{\mu}y_{\mu'}=0$, hence
$X_{\nu,\ft}\ast (w_{\mu'}x_{\mu}w_{\mu}y_{\mu'})=0$ for any $\ft\in\Std(\nu)$. It follows that $\CM_{\Q}^{\lam}(w_{\mu'}x_{\mu}w_{\mu}y_{\mu'})=0$.
Therefore, $S_{\Q}^{\mu}(w_{\mu'}x_{\mu}w_{\mu}y_{\mu'})=0$. On the other hand, since $S_{\Q}^{\mu}\cong x_{\mu}w_{\mu'}y_{\mu'}\Q[\BS_{2n}]$, and by \cite[Lemma 5.7]{KM}, $$
x_{\mu}w_{\mu'}y_{\mu'}(w_{\mu'}x_{\mu}w_{\mu}y_{\mu'})=\Bigl(\prod_{(i,j)\in[\mu]}h_{i,j}^{\mu}\Bigr)
x_{\mu}w_{\mu'}y_{\mu'}\neq 0,
$$
where $h_{i,j}^{\mu}$ is the $(i,j)$-hook length in $[\mu]$, we get
a contradiction. Therefore, $\rho_{\mu}^{\lam}\neq 0$ must imply
that $\mu\unrhd\lam$. Now counting the dimensions, we deduce that
$\CM_{\Q}^{\lam}\cong\oplus_{\lam\unlhd\mu\in 2\PP_n}S_{\Q}^{\mu}$.
In a similar way, we can prove that
$\CM_{\Q}^{\rhd\lam}\cong\oplus_{\lam\lhd\mu\in
2\PP_n}S_{\Q}^{\mu}$. It follows that
$\CM_{\Q}^{\lam}/\CM_{\Q}^{\rhd\lam}\cong S_{\Q}^{\lam}$.\smallskip

We now consider the natural map from $x_{\lam}\Z[\BS_{2n}]$ onto $\CM^{\lam}/\CM^{\rhd\lam}\cong S^{\lam}$.
Since $\Q\otimes_{\Z}M_{\fS_0,\rhd}^{\lam}$ does not contain $S_{\Q}^{\lam}$ as a composition factor, it follows that
(by Proposition \ref{prop210}) the image of $M_{\fS_0,\rhd}^{\lam}$ must be contained in $\CM^{\rhd\lam}$. Therefore
we get a surjective map from $S^{\lam}$ onto $\CM^{\lam}/\CM^{\rhd\lam}\cong S^{\lam}$. This map sends the standard
basis of $S^{\lam}$ to the canonical basis of $\CM^{\lam}/\CM^{\rhd\lam}$. So it must be injective as well, as required.
\end{proof}
\bigskip
\bigskip\bigskip

\section{The $n$-tensor space $V^{\otimes n}$}
\medskip

In this section, we shall use the results obtained in Section 2 and
in \cite{DDH} to give an explicit and characteristic free
description of the annihilator of the $n$-tensor space $V^{\otimes
n}$ in the Brauer algebra $\bb_n(-2m)$. \smallskip

Let $K$ be an arbitrary infinite field. Let $m, n\in\mathbb{N}$.
Let $V$ be a $2m$-dimensional symplectic $K$-vector space. Let
$Sp(V)$ be the corresponding symplectic group, acting naturally on
$V$, and hence on the $n$-tensor space $V^{\otimes n}$ from the
left-hand side. As we mentioned in the introduction, this left
action on $V^{\otimes n}$ is centralized by the specialized Brauer
algebra $\bb_n(-2m)_K:=K\otimes_{\Z}\bb_n(-2m)$, where $K$ is
regarded as $\Z$-algebra by sending each integer $a$ to $a\cdot
1_{K}$. The Brauer algebra $\bb_n(-2m)_K$ acts on the $n$-tensor
space $V^{\otimes n}$ from the right-hand side. Let $\varphi$ be
the natural $K$-algebra homomorphism $$
\varphi:(\bb_n(-2m)_K)^{op}\rightarrow\End_{K}\bigl(V^{\otimes
n}\bigr).
$$

The following is one of the main results in \cite{DDH}, which
generalize earlier results in \cite{B}, \cite{B1}, \cite{B2} for the
case when $K=\mathbb{C}$.

\begin{lem} {\rm (\cite[(1.2)]{DDH})} \label{lm31} Let $K$ be an arbitrary infinite field.
Then $$ \varphi\bigl(\bb_n(-2m)_K\bigr)=\End_{KSp(V)}\bigl(V^{\otimes n}\bigr),$$ and if $m\geq
n$, then $\varphi$ is also injective, and hence an isomorphism.
\end{lem}

By the discussion in \cite[Section 3]{DDH}, $V^{\otimes n}$ is a
tilting module over $KSp(V)$. By \cite[(4.4)]{DPS}, the dimension of
$\End_{KSp(V)}\bigl(V^{\otimes n}\bigr)$ is independent of the
choice of the infinite field $K$. Therefore, the dimension of
$\Ker\varphi:=\bigl\{y\in \bb_n(-2m)_K\bigm|\varphi(y)=0\bigr\}$ is
also independent of the choice of the infinite field $K$. That is,
the dimension of the annihilator of the $n$-tensor space $V^{\otimes
n}$ in the Brauer algebra $\bb_n(-2m)_K$ is independent of the
choice of the infinite field $K$.

\begin{lem} \label{lm32} With the notations as above, we have that $$
\dim(\Ker\varphi)=\sum_{\substack{\lam\in2\PP_n\\ \lam_1>2m}}\dim S^{\lam}.
$$
\end{lem}

\begin{proof} By Lemma \ref{lm22} and Lemma \ref{lm31}, it suffices to consider the case where $K=\mathbb{C}$ and to show that $$
\dim\End_{\mathbb{C} Sp_{2m}(V)}\bigl(V^{\otimes n}\bigr)=\sum_{\substack{\lam\in2\PP_n\\ \lam_1\leq 2m}}\dim S^{\lam}.
$$
Note that $\dim S^{\lam}=\dim S^{\lam'}$, and $$
\End_{\mathbb{C} Sp_{2m}(V)}\bigl(V^{\otimes n}\bigr)\cong \Bigl(\bigl(V^{\otimes n}\bigr)\otimes\bigl(V^{\otimes n}\bigr)^{\ast}\Bigr)^{Sp(V)}\cong \bigl(V^{\otimes 2n}\bigr)^{Sp(V)}.
$$
Therefore, it suffices to show that $$
\dim\bigl(V^{\otimes 2n}\bigr)^{Sp(V)}=\sum_{\substack{\lam\in2\PP_n\\ \lam_1\leq 2m}}\dim S^{\lam'}.
$$
By the well-known Schur-Weyl duality between the general linear group $GL(V)$ and the symmetric group $\BS_{2n}$ on the $2n$-tensor space $V^{\otimes 2n}$, we know that there is a $(GL(V), \BS_{2n})$-bimodules decomposition $$
V^{\otimes 2n}\cong\bigoplus_{\substack{\lam\vdash 2n\\
\ell(\lam)\leq 2m}}\widetilde{\Delta}_{\lam}\otimes S^{\lam},
$$
where $\widetilde{\Delta}_{\lam}$ denotes the irreducible Weyl module with highest weight $\lam$ over $\GL(V)$. Here
we identify $\lam$ with $\lam_1\varepsilon_1+\cdots+\lam_{2n}\varepsilon_{2n}$, $\varepsilon_1, \cdots, \varepsilon_{2n}$ are the fundamental dominant weights of $GL(V)$. It follows that $$
\dim\bigl(V^{\otimes 2n}\bigr)^{Sp(V)}=\sum_{\substack{\lam\vdash 2n\\
\ell(\lam)\leq 2m}}\dim\Bigl((\widetilde{\Delta}_{\lam}\downarrow_{Sp(V)})^{Sp(V)}\Bigr) \dim(S^{\lam}).
$$
By the branching law (see \cite[Proposition 2.5.1]{KT}) from $GL(V)$ to $Sp(V)$, we know that $$
\dim\Bigl((\widetilde{\Delta}_{\lam}\downarrow_{Sp(V)})^{Sp(V)}\Bigr)=1
$$
if $\lam'\in 2\PP_n$; and $0$ otherwise. This proves that $$
\dim\bigl(V^{\otimes 2n}\bigr)^{Sp(V)}=\sum_{\substack{\lam\in2\PP_n\\ \lam_1\leq 2m}}\dim S^{\lam'},
$$
as required.
\end{proof}

Let $a,b$ be two integers such that $0\leq a,b\leq n$ and $a+b$ is
even. Let
$$
I_{a}^{\text{odd}}:=\{1,3,5,\cdots,2a-1\},\quad
I_{b}^{\text{even}}:=\{2,4,6,\cdots,2b\}.
$$
Let $k:=\max\bigl\{2a,2b\bigr\}$. If $k=2a$, we define $\BD_n(a,b)$
to be the set of all the Brauer $n$-diagrams $D$ such that for each
integer $s\in\{1,2,\cdots,2b,2b+1,2b+3,\cdots,2a-1\}$, $D$ connects
the vertex labeled by $s$ with the vertex labeled by $t$ for some
integer $t\in\{1,2,\cdots,2b,2b+1,2b+3,\cdots,2a-1\}\setminus\{s\}$;
and for each integer $s$ with $a+1\leq s\leq n$, $D$ connects the
vertex labeled by $2s-1$ with the vertex labeled by $2s$; and for
each integer $s$ with $1\leq s\leq (a-b)/2$, $D$ connects the vertex
labeled by $2b+4s-2$ with the vertex labeled by $2b+4s$. If $k=2b$,
we define $\BD_n(a,b)$ to be the set of all the Brauer $n$-diagrams
$D$ such that for each integer
$s\in\{1,2,\cdots,2a,2a+2,2a+4,\cdots,2b\}$, $D$ connects the vertex
labeled by $s$ with the vertex labeled by $t$ for some integer
$t\in\{1,2,\cdots,2a,2a+2,2a+4,\cdots,2b\}\setminus\{s\}$; and for
each integer $b+1\leq s\leq n$, $D$ connects the vertex labeled by
$2s-1$ with the vertex labeled by $2s$; and for each integer $1\leq
s\leq (b-a)/2$, $d$ connects the vertex labeled by $2a+4s-3$ with
the vertex labeled by $2a+4s-1$.

\begin{lem} \label{lm33} Let $a,b$ be two integers such that $0\leq a,b\leq n$ and $a+b$ is
even. Suppose that $a+b\geq 2m+2$, then
$$ \sum_{D\in\BD_n(a,b)}D\in\Ker\varphi.
$$
\end{lem}

The proof of Lemma \ref{lm33} is somewhat complicated and will be
given in Section 4.\smallskip

Given any two subsets $A^{(1)}\subseteq\{1,3,\cdots,2n-1\},
A^{(2)}\subseteq\{2,4,\cdots,2n\}$ with $|A^{(1)}|+|A^{(2)}|$ is
even, we set $2n_0=|A^{(1)}|+|A^{(2)}|$, and $$
\bigl\{a_1,a_2,\cdots,a_{2n-2n_0}\bigr\}:=\bigl\{1,2,\cdots,2n\bigr\}\setminus\bigl(A^{(1)}\cup
A^{(2)}\bigr).
$$
Let $(i_1,j_1,i_2,j_2,\cdots,i_{n-n_0},j_{n-n_0})$ be a fixed
permutation of $$\bigl\{a_1,a_2,\cdots,a_{2n-2n_0}\bigr\}.$$ Let
$${\bf i}:=(i_1,i_2,\cdots,i_{n-n_0}), \quad {\bf
j}:=(j_1,j_2,\cdots,j_{n-n_0}).$$ We define $\BD_n^{\bf{i},
\bf{j}}(A^{(1)},A^{(2)})$ to be the set of all the Brauer
$n$-diagrams $D$ such that for each integer $s\in A^{(1)}\cup
A^{(2)}$, $D$ connects the vertex labeled by $s$ with a vertex
labeled by $t$ for some integer $t\in\bigl(A^{(1)}\cup
A^{(2)}\bigr)\setminus \{s\}$, and for each integer $s$ with $1\leq
s\leq n-n_0$, $D$ connects the vertex labeled by $i_s$ with the
vertex labeled by $j_s$. Note that the set $\BD_n(a,b)$ we defined
before is a special case of the set
$\BD_n^{\bf{i},\bf{j}}(A^{(1)},A^{(2)})$ we defined here.

\addtocounter{cor}{3}
\begin{cor} \label{cor34} With the notations as above and suppose that
$2n_0=|A^{(1)}|+|A^{(2)}|\geq 2m+2$, then we have
$$ \sum_{D\in\BD_n^{\bf{i},\bf{j}}(A^{(1)},A^{(2)})}D\in\Ker\varphi.
$$
\end{cor}

\begin{proof} Let $n_1=|A^{(1)}|, n_2=|A^{(2)}|$. If
$n_1\geq n_2$, then for any Brauer diagram
$D\in\BD_n^{\bf{i},\bf{j}}(A^{(1)},A^{(2)})$, there exist at least
$\frac{n_1-n_2}{2}$ bottom horizontal edges between the vertices
labeled by the integers in the following set $$
\bigl\{2,4,6,\cdots,2n\bigr\}\setminus A^{(2)}.
$$
As a result, we deduce that there exist elements
$\sigma_{A^{(1)}}\in\BS_{(1,3,\cdots,2n-1)}$,
$\sigma_{A^{(2)}}\in\BS_{(2,4,\cdots,2n)}$ and a Brauer diagram
$D_1$, such that

\begin{enumerate}
\item[(1)] for any integer $a$ with $1\leq a\leq n_1$, $D_1$ connects the
vertex labeled by $2a-1$ with the vertex labeled by $2a$.
\item[(2)] $$
\sigma_{A^{(1)}}\Bigl(\sum_{D\in\BD_n^{\bf{i},\bf{j}}(A^{(1)},A^{(2)})}D\Bigr)\sigma_{A^{(2)}}=
D_1\Bigl(\sum_{D\in\BD_n(|A^{(1)}|,|A^{(2)}|)}D\Bigr).
$$
\end{enumerate}
In this case, since both $\varphi(\sigma_{A^{(1)}})$ and
$\varphi(\sigma_{A^{(2)}})$ are invertible, it follows directly from
Lemma \ref{lm33} that $
\sum_{D\in\BD_n^{\bf{i},\bf{j}}(A^{(1)},A^{(2)})}D\in\Ker\varphi$.

If $n_1\leq n_2$, then for any Brauer diagram
$D\in\BD_n^{\bf{i},\bf{j}}(A^{(1)},A^{(2)})$, there exist at least
$\frac{n_2-n_1}{2}$ top horizontal edges between the vertices
labeled by the integers in the following set $$
\bigl\{1,3,5,\cdots,2n-1\bigr\}\setminus A^{(1)}.
$$
As a result, we deduce that there exist elements
$\sigma_{A^{(1)}}\in\BS_{(1,3,\cdots,2n-1)}$,
$\sigma_{A^{(2)}}\in\BS_{(2,4,\cdots,2n)}$ and a Brauer diagram
$D_2$, such that

\begin{enumerate}
\item[(3)] for any integer $a$ with $1\leq a\leq n_2$, $D_2$ connects the
vertex labeled by $2a-1$ with the vertex labeled by $2a$.
\item[(4)] $$
\sigma_{A^{(1)}}\Bigl(\sum_{D\in\BD_n^{\bf{i},\bf{j}}(A^{(1)},A^{(2)})}D\Bigr)\sigma_{A^{(2)}}=
D_2\Bigl(\sum_{D\in\BD_n(|A^{(1)}|,|A^{(2)}|)}D\Bigr).
$$
\end{enumerate}
By the same argument as before, we deduce that $
\sum_{D\in\BD_n^{\bf{i},\bf{j}}(A^{(1)},A^{(2)})}D\in\Ker\varphi$ in
this case. This completes the proof of the corollary.
\end{proof}

The following is the main result of this section, which gives an
explicit and characteristic free description of the annihilator of
the $n$-tensor space $V^{\otimes n}$ in the Brauer algebra
$\bb_n(-2m)$.

\addtocounter{thm}{4}
\begin{thm} \label{thm35} With the notations as in Lemma \ref{lm31} and Lemma \ref{lm32}, we have that $$
\Ker\varphi=\CM_{K}^{(2m+2,2^{n-m-1})},
$$
where
$(2m+2,2^{n-m-1}):=(2m+2,\underbrace{2,\cdots,2}_{\text{$n-m-1$
copies}})$, $\CM_{K}^{(2m+2,2^{n-m-1})}$ is the right
$K[\BS_{2n}]$-module associated to $(2m+2,2^{n-m-1})$ as defined in
Section 2. In particular, $\Ker\varphi$ is a $\BS_{2n}$-submodule.
\end{thm}

\begin{proof} It is easy to see that for any partition $\mu\in 2\PP_n$, $\mu\unrhd (2m+2,2^{n-m-1})$ if and only if
$\mu_1>2m$. Therefore, $$
\dim\CM_{K}^{(2m+2,2^{n-m-1})}=\sum_{\substack{\lam\in2\PP_n\\ \lam_1>2m}}\dim S^{\lam}.
$$
Applying Lemma \ref{lm31} and Lemma \ref{lm32}, we see that to prove
this theorem, it suffices to show that
$\CM_{K}^{(2m+2,2^{n-m-1})}\subseteq\Ker\varphi$. Equivalently, it
suffices to show that for any partition
$\lam=(\lam_1,\cdots,\lam_s)\in 2\PP_n$ satisfying $\lam_1>2m$, and
any $w\in\BS_{2n}$, $\varphi(X_{\lam}\ast w)=0$.\smallskip

By the definition of the element $X_{\lam}$, the action ``$\ast$"
and the multiplication rule of Brauer diagrams, we deduce that
$$
X_{\lam}\ast
w=\sum_{\bf{i},\bf{j}}\sum_{D\in\BD_n^{\bf{i},\bf{j}}(A^{(1)},A^{(2)})}D,
$$
where $$\begin{aligned}
A^{(1)}&:=\Bigl\{(i)w\Bigm|i=1,2,3,\cdots,\lam_1\Bigr\}\bigcap\Bigl\{1,3,5,\cdots,2n-1\Bigr\},\\
A^{(2)}&:=\Bigl\{(i)w\Bigm|i=1,2,3,\cdots,\lam_1\Bigr\}\bigcap\Bigl\{2,4,6,\cdots,2n\Bigr\},
\end{aligned}
$$
and $|A^{(1)}|+|A^{(2)}|=2n_0=\lambda_1$, and ${\bf
i}:=(i_1,i_2,\cdots,i_{n-n_0})$, ${\bf
j}:=(j_1,j_2,\cdots,j_{n-n_0})$ such that
$(i_1,j_1,i_2,j_2,\cdots,i_{n-n_0},j_{n-n_0})$ is a permutation of
the integers in $\bigl\{1,2,\cdots,2n\bigr\}\setminus(A^{(1)}\cup
A^{(2)})$.\smallskip

We now apply Corollary \ref{cor34}. It follows immediately that
$\varphi(X_{\lam}\ast w)=0$ as required. This completes the proof of
the theorem.
\end{proof}

\addtocounter{rem}{5}
\begin{rem}\label{rem36} Let $V_{\Z}$ be a free $\Z$-module with
basis $\{v_1,v_2,\cdots,v_{2m}\}$. For any $\Z$-algebra $R$, we
define $V_{R}:=R\otimes_{\Z}V_{\Z}$. The same formulae (see Section
1) define an action of the algebra $\bb_n(-2m)$ on $V_{\Z}^{\otimes
n}$, and hence an action of $\bb_n(-2m)_{R}$ on $V_{R}^{\otimes n}$.
Let $S_{R}^{sy}(m,n)$ (see \cite[Section 2]{DDH} and \cite{Oe}) be
the symplectic Schur algebra over $R$. If $R$ is a field, then
$S_{R}^{sy}(m,n)$ is a quasi-hereditary algebras over $R$, and
$V_{R}^{\otimes n}$ is a tilting module over $S_{R}^{sy}(m,n)$.
Applying \cite[(4.4)]{DPS}, we know that, for any $\Z$-algebra $R$,
there is a canonical isomorphism
$$ R\otimes_{\Z}\End_{S_{\Z}^{sy}(m,n)}\Bigl(V_{\Z}^{\otimes
n}\Bigr)\cong \End_{S_{R}^{sy}(m,n)}\Bigl(V_{R}^{\otimes n}\Bigr).
$$
Note that $$
\varphi\bigl(\bb_n(-2m)_R\bigr)\subseteq\End_{S_{R}^{sy}(m,n)}\Bigl(V_{R}^{\otimes
n}\Bigr)
$$
By the main result in \cite{DDH}, we know that the above inclusion
``$\subseteq$" can be replaced by ``=" when $R=K$ is an infinite
field $K$. In fact, this is always true for any $\Z$-algebra $R$ (by
using some localization argument in commutative algebras). As a
consequence, our Theorem \ref{thm35} is also always true if we
replace the infinite field $K$ by any $\Z$-algebra $R$.
\end{rem}

\bigskip\bigskip
\section{Proof of Lemma 3.3}
\medskip

We shall first fix some notations and convention. Note that the
element $\sum_{D\in\BD_n(a,b)}D$ in Lemma \ref{lm33} actually lies
in $\bb_n$, we can choose to work inside the Brauer algebra
$\bb_n(-2m)_{\mathbb{C}}$ in this section. Furthermore, throughout
this section, we shall recover our original way of labeling of
vertices in each Brauer $n$-diagram. That is, the vertices in each
row of a Brauer $n$-diagram will be labeled by the indices
$1,2,\cdots,n$ from left to right. This way of labeling is more
convenient when we need to express each Brauer diagram in terms of
the standard generators $s_i, e_i$ for $1\leq i\leq n-1$ and to
consider the action of Brauer diagrams on the $n$-tensor space
$V^{\otimes n}$.\smallskip

Let $f$ be an integer with $0\leq f\leq [n/2]$, where $[n/2]$ is
the largest non-negative integer not bigger than $n/2$. Let
$\nu=\nu_{f}:=((2^f), (n-2f))$, where
$(2^f):=(\underbrace{2,2,\cdots,2}_{\text{$f$ copies}})$ and
$(n-2f)$ are considered as partitions of $2f$ and $n-2f$
respectively. In general, a bipartition of $n$ is a pair
$(\lam^{(1)},\lam^{(2)})$ of partitions of numbers $n_1$ and $n_2$
with $n_1+n_2=n$. The notions of Young diagram, bitableaux, etc.,
carry over easily. Let $\ft^{\nu}$ be the standard $\nu$-bitableau
in which the numbers $1,2,\cdots,n$ appear in order along
successive rows of the first component tableau, and then in order
along successive rows of the second component tableau. We define
$$
\D_{\nu}:=\Biggl\{d\in\BS_n\Biggm|\begin{matrix}\text{$(\ft^{(1)},\ft^{(2)})=\ft^{\nu}d$
is row standard and the first}\\
\text{ column of $\ft^{(1)}$ is an increasing sequence when}\\
\text{read from top to bottom}\\
\end{matrix}\Biggr\}.
$$
For each partition $\lambda$ of $n-2f$, we denote by
$\Std_{2f}(\lambda)$ the set of all the standard
$\lambda$-tableaux with entries in $\{2f+1,\cdots,n\}$. The
initial tableau $\ft^{\lam}$ in this case has the numbers
$2f+1,\cdots,n$ in order along successive rows. Again, for each
$\ft\in\Std_{2f}(\lambda)$, let $d(\ft)$ be the unique element in
$\BS_{\{2f+1,\cdots,n\}}\subseteq\BS_n$ with
$\ft^{\lam}d(\ft)=\ft$. Let $\sig\in\BS_{\{2f+1,\cdots,n\}}$ and
$d_1,d_2\in\mathcal{D}_{\nu}$. Then $d_1^{-1}e_1e_3\cdots
e_{2f-1}\sig d_2$ corresponds to the Brauer $n$-diagram where the
top horizontal edges connect $(2i-1)d_1$ and $(2i)d_1$, the bottom
horizontal edges connect $(2i-1)d_2$ and $(2i)d_2$, for
$i=1,2,\cdots,f$, and the vertical edges connects $(j)d_1$ with $
(j)d_2$ for $j=2f+1,2f+2,\cdots,n$.

\begin{lem} {\rm (\cite[Corollary 3.3]{DDH})} With the above notations, the set
$$ \biggl\{d_1^{-1}e_1e_3\cdots e_{2f-1}\sig
d_2\biggm|\begin{matrix}\text{$0\leq f\leq [n/2]$,
$\sig\in\BS_{\{2f+1,\cdots,n\}}$, $d_1,
d_2\in\D_{\nu}$,}\\ \text{where $\nu:=((2^f), (n-2f))$}\\
\end{matrix}\biggr\}.$$
is a basis of Brauer algebra $B_{n}(x)_R$, which coincides with
the natural basis given by Brauer $n$-diagrams.
\end{lem}
Given an element $d_1^{-1}e_1e_3\cdots e_{2f-1}\sig d_2$ as above,
let $D$ be its representing Brauer $n$-diagram. Let
$v_{\ui}:=v_{i_1}\otimes v_{i_2}\otimes\cdots\otimes v_{i_n}$ be a
simple $n$-tensor in $V^{\otimes n}$.

\begin{lem} \label{lm42} With the notations as above, we have that
$$v_{\ui}D=(-1)^{\ell(d_1^{-1}\sigma d_2)}\bigl(v_{\ui}\circ
D)\bigr),$$ where $v_{\ui}\circ D$ can be described as
follows:\smallskip

1) If $(j)d_1^{-1}\sigma d_2=(k)$ for
$j\in\bigl\{(2f+1)d_1,(2f+2)d_1,\cdots,(n)d_1\bigr\}$, then the
$k$th position of $v_{\ui}\circ D$ is $v_{i_{j}}$;\smallskip

2) For each $1\leq j\leq f$, the $((2j-1)d_2, (2j)d_2)$th position
of $v_{\ui}\circ D$ is the following sum:
$$ \epsilon_{i_{(2j-1)d_1}, i_{(2j)d_1}} \sum_{k=1}^{m}(v_{k'}\otimes v_{k}-
v_{k}\otimes v_{k'}).
$$
\end{lem}

\addtocounter{rem}{2}
\begin{rem} \label{rem41} Intuitively, the action of the Brauer $n$-diagram $D$ on
$v_{\ui}$ can be thought as follows. Let
$(a_1,b_1),\cdots,(a_f,b_f)$ be the set of all the horizontal
edges in the top row of $D$, where $a_s<b_s$ for each $s$ and
$a_1<a_2<\cdots<a_f$. Let $(c_1,d_1),\cdots,(c_f,d_f)$ be the set
of all the horizontal edges in the bottom row of $D$, where
$c_s<d_s$ for each $s$ and $c_1<c_2<\cdots<c_f$. Then for each
$1\leq j\leq f$, the $(c_j, d_j)$th position of $v_{\ui}\circ D$
is the following sum:
$$ \epsilon_{i_{a_j}, i_{b_j}} \sum_{k=1}^{m}(v_{k'}\otimes v_{k}-
v_{k}\otimes v_{k'}).
$$
We list those vertices in the top row of $D$ which are not connected
with horizontal edges from left to right as
$i_{k_{2f+1}},i_{k_{2f+2}},\cdots,i_{k_{n}}$. Then, for each integer
$j$ with $2f+1\leq j\leq n$, the $(j\sigma d_2)$th position of
$v_{\ui}\circ D$ is $v_{i_{k_{j}}}$.
\end{rem}

We fix an arbitrary element $d_2\in\D_{\nu}$. We define $$
\BD^{(f)}(n;d_2):=\bigl\{d_1^{-1}e_1e_3\cdots e_{2f-1}\sigma
d_2\bigm|\text{$d_1\in\D_{\nu},
\sigma\in\BS_{\{2f+1,2f+2,\cdots,n\}}$}\bigr\}.
$$
Note that $\BD^{(f)}(n;d_2)$ consists of all the Brauer
$n$-diagrams whose bottom horizontal edges are $$
((1)d_2,(2)d_2),((3)d_2,(4)d_2),\cdots,((2f-1)d_2,(2f)d_2).
$$

\addtocounter{lem}{1}
\begin{lem} \label{lm43} Let $f$ be an integer with $0\leq f\leq [n/2]$. Let $d_2\in\D_f$.
Then for any $\sigma\in\BS_{n}$, $$
\sigma\Bigl(\sum_{D\in\BD^{(f)}(n;d_2)}D\Bigr)=\sum_{D\in\BD^{(f)}(n;d_2)}D.
$$
\end{lem}

\begin{proof} It suffices to show that for each integer $1\leq i<n$,
\addtocounter{equation}{4}
\begin{equation}\label{equa41}
s_i\Bigl(\sum_{D\in\BD^{(f)}(n;d_2)}D\Bigr)=\sum_{D\in\BD^{(f)}(n;d_2)}D.
\end{equation}
In fact, for $D, D'\in\BD^{(f)}(n;d_2)$ with $D\neq D'$, it is
clear that $s_iD\neq s_iD'$, and both $s_iD$ and $s_iD'$ are still
lie in $\BD^{(f)}(n;d_2)$. Now counting the number of Brauer
$n$-diagrams occurring in both sides, we prove (\ref{equa41}) and
hence also prove the lemma.
\end{proof}

Similarly, we fix an arbitrary element $d_1\in\D_{\nu}$ and define
$$ \BD^{(f)}(d_1;n):=\bigl\{d_1^{-1}e_1e_3\cdots e_{2f-1}\sigma
d_2\bigm|\text{$d_2\in\D_{\nu},
\sigma\in\BS_{\{2f+1,2f+2,\cdots,n\}}$}\bigr\}.
$$
Then $\BD^{(f)}(d_1;n)$ consists of all the Brauer $n$-diagrams
whose top horizontal edges are $$
((1)d_1,(2)d_1),((3)d_1,(4)d_1),\cdots,((2f-1)d_1,(2f)d_1).
$$
The following result can be proved in the same way as Lemma
\ref{lm43}.

\addtocounter{lem}{1}
\begin{lem} \label{lm44} Let $f$ be an integer with $0\leq f\leq [n/2]$. Let $d_1\in\D_f$.
Then for any $\sigma\in\BS_{n}$, $$
\Bigl(\sum_{D\in\BD^{(f)}(d_1;n)}D\Bigr)\sigma=\sum_{D\in\BD^{(f)}(d_1;n)}D.
$$
\end{lem}

Let $\ui=(i_1,i_2,\cdots,i_n)$, where $1\leq i_j\leq 2m$ for each
$j$. An ordered pair $(s,t)$ ($1\leq s<t\leq n$) is called a {\it
symplectic pair} in $\ui$ if $i_s=i'_t$. Two ordered pairs $(s,t)$
and $(u,v)$ are called disjoint if
$\bigl\{s,t\bigr\}\cap\bigl\{u,v\bigr\}=\emptyset$. We define the
{\it symplectic length} $\ell_s(v_{\ui})$ to be the maximal number
of disjoint symplectic pairs $(s,t)$ in $\ui$. We now consider a
special case of Lemma \ref{lm33}.

\addtocounter{prop}{6}
\begin{prop} \label{prop41} We have that $$
\sum_{D\in\BD_n(m+1,m+1)}D\in\Ker\varphi.
$$
\end{prop}

\begin{proof} By the above discussion and the definition of
$\BD_n(m+1,m+1)$, any Brauer diagram $d\in \BD_n(m+1,m+1)$ only
acts on the first $m+1$ components of any simple $n$-tensor
$v_{i_1}\otimes v_{i_2}\otimes\cdots\otimes v_{i_n}\in V^{\otimes
n}$. Therefore, to show that
$\sum_{D\in\BD_n(m+1,m+1)}D\in\Ker\varphi$, we can assume without
loss of generality that $n=m+1$. Note that $$\begin{aligned}
\sum_{D\in\BD_n(m+1,m+1)}D&=\sum_{0\leq f\leq
[n/2]}\sum_{d_{2}\in\D_{{\nu\!_{f}}}}\sum_{D\in\BD^{(f)}(n;d_2)}D\\
&=\sum_{0\leq f\leq
[n/2]}\sum_{d_{1}\in\D_{{\nu\!_{f}}}}\sum_{D\in\BD^{(f)}(d_1;n)}D.\end{aligned}
$$
Suppose that $\sum_{D\in\BD_n(m+1,m+1)}D\not\in\Ker\varphi$. Then
there exists a simple $n$-tensor $v_{\ui}=v_{i_1}\otimes
v_{i_2}\otimes\cdots\otimes v_{i_n}\in V^{\otimes n}$, such that
$$ v_{\ui}\sum_{D\in\BD_n(m+1,m+1)}D\neq
0.$$ Let $\ui:=(i_1,\cdots,i_n)$. Suppose that $\ell_{s}(\ui)=f$
for some integer $0\leq f\leq [n/2]$.

By Lemma \ref{lm43}, we have that for any $\sigma\in\BS_n$, $$
\sigma\sum_{D\in\BD_n(m+1,m+1)}D=\sum_{D\in\BD_n(m+1,m+1)}D.
$$
Therefore, we can assume without loss of generality that
$i_{2s-1}=i'_{2s}<i_{2s}$ for each integer $1\leq s\leq f$,
$i_1\leq i_3\leq \cdots\leq i_{2f-1}$ and $i_{2f+1}\leq
i_{2f+2}\leq\cdots\leq i_n$. Furthermore, if $i_j=i_k$ for some
integers $j\neq k$, then $$\begin{aligned}
v_{\ui}\sum_{D\in\BD_n(m+1,m+1)}D&=
v_{\ui}\frac{1+s_{(j,k)}}{2}\sum_{D\in\BD_n(m+1,m+1)}D=0,
\end{aligned}$$
where $s_{(j,k)}$ denotes the transposition $(j,k)$ in $\BS_n$, and
we have used the fact that the length of $s_{(j,k)}$ is an odd
integer. Therefore, we can deduce that $i_1,i_2,\cdots,i_n$ are
pairwise distinct. Hence, $i_{2s-1}=i'_{2s}<i_{2s}$ for each integer
$1\leq s\leq f$, $i_1<i_3<\cdots<i_{2f-1}$ and
$i_{2f+1}<i_{2f+2}<\cdots<i_n$.

By Lemma \ref{lm42}, we get that \addtocounter{equation}{2}
\begin{equation}\label{equa44} v_{\ui}\sum_{0\leq g\leq
f}\sum_{d_{1}\in\D_{\nu\!_g}}\sum_{D\in\BD^{(g)}(d_1;n)}D\neq
0.\end{equation}

By Lemma \ref{lm44}, it is easy to see that for any $w\in\BS_n$,
$$
v_{\ui}\sum_{0\leq g\leq
f}\sum_{d_{1}\in\D_{\nu\!_g}}\sum_{D\in\BD^{(g)}(d_1;n)}Dw=v_{\ui}\sum_{0\leq
g\leq f}\sum_{d_{1}\in\D_{\nu\!_g}}\sum_{D\in\BD^{(g)}(d_1;n)}D.
$$
In particular, by the same argument as before, any simple $n$-tensor
$v_{b_1}\otimes\cdots\otimes v_{b_n}$ with $b_j=b_k$ for some
integers $j\neq k$ can not appear with nonzero coefficient in the
expansion of the left-hand side of (\ref{equa44}). Therefore if a
simple $n$-tensor $v_{\ub}$ appears with nonzero coefficient in the
expansion of the left-hand side of (\ref{equa44}), then
$b_1,\cdots,b_n$ must be pairwise distinct. Applying Lemma
\ref{lm43} and Lemma \ref{lm44} again, we can further assume
(without loss of generality) that there exists such a $v_{\ub}$
which appears with nonzero coefficient in the expansion of the
left-hand side of (\ref{equa44}), such that $b_1,\cdots,b_n$ are
pairwise distinct, $\ell_s(\ub)=f$, $b_{2s-1}=b'_{2s}<b_{2s}$ for
each integer $1\leq s\leq f$, $b_1<b_3<\cdots<b_{2f-1}$, and there
exists integer $0\leq r\leq f$ such that $b_{t}=i_{t}$ for each
integer $2r+1\leq t\leq n$, and
$\{b_1,b_2,\cdots,b_{2r}\}\cap\{i_1,i_2,\cdots,i_{2r}\}=\emptyset$.
\smallskip

Let $g$ be an integer with $0\leq g\leq f$, $d_1\in\D_{\nu\!_g}$,
$D\in\BD^{(g)}(d_1;n)$, where $$ D=d_1^{-1}e_1e_3\cdots
e_{2g-1}\sigma d_2,\,\,\sigma\in\BS_{\{2g+1,2g+2,\cdots,n\}},\,\,
d_2\in\D_{\nu_g}.
$$
We claim that $v_{\ub}$ appears with nonzero coefficient in the
expansion of $v_{\ui}D$ if and only if \smallskip

1) $g\geq r$, $\sigma=1$, and\medskip

2) the horizontal edges in the top row of $D$ are of the form
$$\begin{aligned}
&(1,2),(3,4),\cdots,(2r-1,2r),(2a_1-1,2a_1),(2a_2-1,2a_2),\cdots,\\
&\qquad\qquad\qquad(2a_{g-r}-1,2a_{g-r}),\end{aligned}
$$
where $a_1,\cdots,a_{g-r}$ are some integers satisfying $r+1\leq
a_1<a_2<\cdots<a_{g-r}\leq f$, and
\smallskip

3) the horizontal edges in the bottom row of $D$ is the same as
those in the top row of $D$, i.e., $d_2=d_1$.\smallskip

In fact, for any Brauer diagram $D$ satisfying the above conditions
1), 2), 3), by Remark \ref{rem41}, $v_{\ub}$ does appear with
nonzero coefficient in the expansion of $v_{\ui}D$, the coefficient
is $(-1)^g$; on the other hand, suppose that $v_{\ub}$ appears with
nonzero coefficient in the expansion of $v_{\ui}D$. By our
assumption on $\ui$ and $\ub$, it is easy to see that the tensor
factor $v_{b_{2s-1}}\otimes v_{b_{2s}}$ with $1\leq s\leq r$ can
only be produced through the action of $e_{2t-1}$ for some $1\leq
t\leq g$. This implies that $g\geq r$. For each integer $j$ with
$2g+1\leq j\leq n$, by Remark \ref{rem41}, the action of $D$ on
$v_{\ui}$ move the vector in the $jd_1$th position of $v_{\ui}$
(i.e., $v_{i_{jd_1}}$) to the $(j\sigma d_2)$th position. By our
assumption on $\ui, \ub$ again, we deduce that $jd_1=j\sigma d_2$.
But by the definition of $\D_g$,
$$\begin{aligned} &
(2g+1)d_1<(2g+2)d_1<\cdots<(n)d_1,\\
& (2g+1)d_2<(2g+2)d_2<\cdots<(n)d_2.\end{aligned}
$$
It follows that $\sigma=1$, and $jd_1=jd_2$ for any $2g+1\leq j\leq
n$. Now the remaining statements of our claim follows easily from
the fact that $\sigma=1$, our assumption on $\ui$ and $\ub$ as well
as Remark \ref{rem41}.

\smallskip
Therefore, the coefficient of $v_{\ub}$ in the expansion of
$$\sum_{d_{1}\in\D_{\nu\!_g}}\sum_{D\in\BD^{(g)}(d_1;n)}v_{\ui}D$$
is equal to
$$ (-1)^g\begin{pmatrix}f-r\\
g-r\end{pmatrix}.
$$
Hence the coefficient of $v_{\ub}$ in the left-hand side of
(\ref{equa44}) is $$ \sum_{r\leq g\leq
f}(-1)^g\begin{pmatrix}f-r\\ g-r\end{pmatrix}=0,
$$
a contradiction. This completes the proof of the proposition.
\end{proof}

Next we consider a more general situation than Proposition
\ref{prop41}, which is still a special case of Lemma \ref{lm33}.

\addtocounter{prop}{1}
\begin{prop} \label{prop42} Let $a,b$ be two integers such that $0\leq a,b\leq n$ and $a+b$ is
even. Suppose that $a+b\geq 2m+2$ and $a\geq b$, then
$$ \sum_{D\in\BD_n(a,b)}D\in\Ker\varphi.
$$
\end{prop}

\begin{proof} By the assumption that $a\geq b$ and the definition of $\BD_n(a,b)$, any Brauer diagram
$d\in\BD_n(a,b)$ only acts on the first $a$ components of any
simple $n$-tensor in $V^{\otimes n}$. Therefore, we can assume
without loss of generality that $n=a$. Also, because of
Proposition \ref{prop41}, we can assume that $n=a>b$. In
particular, $n=a\geq m+2, n+b\geq 2m+2$. Suppose that
$\sum_{D\in\BD_n(n,b)}D\not\in\Ker\varphi$. Then there exists a
simple $n$-tensor $v_{\ui}\in V^{\otimes n}$ such that
\addtocounter{equation}{1}
\begin{equation} \label{equa45} v_{\ui}\sum_{D\in\BD_n(n,b)}D\neq 0.
\end{equation}
Applying Lemma \ref{lm43}, we know that for any $\sigma\in\BS_n$,
\begin{equation}\label{equa46}
\sigma\sum_{D\in\BD_n(n,b)}D=\sum_{D\in\BD_n(n,b)}D.
\end{equation}
Therefore, by the same argument as before, we deduce that
$i_1,\cdots,i_n$ are pairwise distinct.\smallskip

Let $f:=(n-b)/2$. We define $$
\Sigma_f=\biggl\{(\fa,\fb):=((a_1,b_1),\cdots,(a_{f},b_{f}))\biggm|
\begin{matrix}\text{$1\leq a_1<\cdots<a_f\leq n$,}\\
\text{$a_i<b_i$ for each $1\leq i\leq f$}
\end{matrix}
\biggr\}.
$$
For each $(\fa,\fb)\in\Sigma_f$, we define $\BD_{\fa,\fb}(n,b)$ to
be the set of all the Brauer $n$-diagrams in $\BD_n(n,b)$ whose
rightmost $f$ horizontal edges in the top row are exactly those
pairs in $(\fa,\fb)$. It is clear that $$
\BD_n(n,b)=\bigsqcup_{(\fa,\fb)\in\Sigma_f}\BD_{\fa,\fb}(n,b).
$$
Therefore,
\begin{equation} \label{equa47} v_{\ui}\sum_{D\in\BD_{\fa,\fb}(n,b)}D\neq 0,
\end{equation}
for some $(\fa,\fb)\in\Sigma_f$. Henceforth, we fix such an
$(\fa,\fb)$. We list the integers in the set $\{k|a_1\leq k\leq
n\}\setminus\{a_i, b_i|1\leq i\leq f\}$ as
$t_1,t_2,\cdots,t_{b+1-a_1}$ such that $t_1<t_2<\cdots<t_{b+1-a_1}$.
We define
$$\fI_{\fa,\fb}:=\Bigl\{\uj:=(j_1,\cdots,j_{b+1-a_1})\Bigm|1\leq j_1<\cdots<j_{b+1-a_1}\leq
b\Bigr\}.$$ For each element $\sigma\in\BS_{b+1-a_1}$, we define
$\BD_{\fa,\fb}^{\sigma,\uj}(n,b)$ to be the set of all the Brauer
$n$-diagrams in $\BD_{\fa,\fb}(n,b)$ whose vertical edges connects
the vertex labeled by $t_s$ in the top row with the vertex labeled
by $j_{s\sigma}$ in the bottom row for each integer $s$ with $1\leq
s\leq b+1-a_1$. It is clear that $$
\BD_{\fa,\fb}(n,b)=\bigsqcup_{\substack{\sigma\in\BS_{b+1-a_1}\\
\uj\in\fI_{\fa,\fb}}}\BD_{\fa,\fb}^{\sigma,\uj}(n,b).
$$

Now let $\uk:=(k_1,\cdots,k_n)$ be such that $v_{\uk}$ appears
with nonzero coefficient in the expansion of (\ref{equa47}). Using
the same argument in the proof of Lemma \ref{lm44}, one can show
that for any $\sigma\in\BS_{b}$, $$
\Bigl(\sum_{D\in\BD_n(n,b)}D\Bigr)\sigma=\sum_{D\in\BD_n(n,b)}D.
$$
As a result, we deduce that $k_1,\cdots,k_b$ are pairwise distinct.
Furthermore, we can assume (if necessary, we replace $v_{\ui}$ by
$v_{\ui}\sigma$ for some $\sigma\in\BS_{n}$) that
$$
k_{b-s+1}=i_{t_{b+2-a_1-s}},\,\,\text{for $s=1,2,\cdots,b+1-a_1$.}
$$
It follows that $v_{\uk}$ appears with nonzero coefficient in the
expansion of \begin{equation}\label{equa48}
v_{\ui}\sum_{D\in\BD_{\fa,\fb}^{(0)}(n,b)}D,
\end{equation}
where
$\BD_{\fa,\fb}^{(0)}:=\BD_{\fa,\fb}^{1,({a_1},{a_1+1},\cdots,{b})}$.
We divide the remaining proof into two cases:
\smallskip

{\it Case 1.} $\ell_s(\ui)=(n-b)/2$. Then it follows from
(\ref{equa48}) and our assumption on $\ui$ that
$\{i_1,\cdots,i_{a_1-1}\}=\{k_1,\cdots,k_{a_1-1}\}$, and the
vectors $v_{i_1},\cdots,v_{i_{a_1-1}}$ lie in a subspace with
dimension
$$\leq m-\frac{n-b}{2}-(b+1-a_1)\leq a_1-2.$$ Hence $i_s=i_t$
for some $1\leq s\neq t\leq a_1-1$. But (\ref{equa48}) implies
that $$ (v_{i_1}\otimes\cdots\otimes
v_{i_{a_1-1}})\sum_{\sigma\in\BS_{a_{1}-1}}\sigma\neq 0,
$$
a contradiction, as required.
\smallskip

{\it Case 2.} $\ell_s(\ui)>(n-b)/2$. Let $W$ be the symplectic
subspace generated by $v_{k_1},\cdots,v_{k_{a_1-1}}$. Then
(\ref{equa48}) implies that $\dim W\leq 2(a_1-2)$. Hence
$\widetilde{m}:=(\dim W)/2\leq a_1-2$. Note that if $i_s\not\in W$
for some $1\leq s\leq a_1-1$, then we must have that $i_s=i'_t$ for
some integer $t$ with $1\leq t\leq a_1-1$. Furthermore, in this
case, if we replace the tensor factors $v_{i_s}, v_{i_t}$ in
$v_{\ui}$ by $v_{k_1}, v_{k'_1}$ respectively, then the coefficient
of $v_{\uk}$ in (\ref{equa48}) changes at most one sign. Therefore,
(\ref{equa48}) implies that there exists a simple $(a_1-1)$-tensor
$\widetilde{v}$ in $W^{\otimes a_1-1}$ such that $v_{k_1}\otimes
\cdots\otimes v_{k_{a_1-1}}$ appears with non-zero coefficient in
$$
\widetilde{v}\sum_{D\in\BD_{a_1-1}(\widetilde{m},\widetilde{m})}D,
$$
where the above element $D$ is understood as element in the Brauer
algebra $\bb_n(-2\widetilde{m})$, acting on the $(a_1-1)$-tensor
space $W^{\otimes a_1-1}$. This is impossible by Proposition
\ref{prop41}. Hence we complete the proof of this
proposition.\end{proof}

Finally, thanks to Proposition \ref{prop42}, to complete the proof
of Lemma \ref{lm33}, we only need to prove the following
proposition.

\addtocounter{prop}{4}
\begin{prop} \label{prop43} Let $a,b$ be two integers such that $0\leq a,b\leq n$ and $a+b$ is
even. Suppose that $a+b\geq 2m+2$ and $b\geq a$, then
$$ \sum_{D\in\BD_n(a,b)}D\in\Ker\varphi.
$$
\end{prop}

\begin{proof} By the assumption that $b\geq a$ and the definition of $\BD_n(a,b)$, any Brauer diagram
$d\in\BD_n(a,b)$ only acts on the first $b$ components of any
simple $n$-tensor in $V^{\otimes n}$. Therefore, we can assume
without loss of generality that $n=b$. Also, because of
Proposition \ref{prop41}, we can assume that $n=b>a$. In
particular, $n=b\geq m+2, n+a\geq 2m+2$. Suppose that
$\sum_{D\in\BD_n(a,n)}D\not\in\Ker\varphi$. Then there exists a
simple $n$-tensor $v_{\ui}\in V^{\otimes n}$ such that
\addtocounter{equation}{1}
\begin{equation} \label{equa49} v_{\ui}\sum_{D\in\BD_n(a,n)}D\neq 0.
\end{equation}
It follows that $i_{a+2s-1}=i'_{a+2s}$ for each integer $s$ with
$1\leq s\leq (n-a)/2$. In particular, $\ell_s(\ui)\geq (n-a)/2$.
Using the same argument in the proof of Lemma \ref{lm44}, one can
show that for any $\sigma\in\BS_{a}$, $$
\sigma\Bigl(\sum_{D\in\BD_n(a,n)}D\Bigr)=\sum_{D\in\BD_n(a,n)}D.
$$
As a result, we deduce that $i_1,\cdots,i_a$ are pairwise
different.\smallskip

Let $f:=(n-a)/2$. We define $$
\Sigma_f=\biggl\{(\fa,\fb):=((a_1,b_1),\cdots,(a_{f},b_{f}))\biggm|
\begin{matrix}\text{$1\leq a_1<\cdots<a_f\leq n$,}\\
\text{$a_i<b_i$ for each $1\leq i\leq f$}
\end{matrix}
\biggr\}.
$$
For each $(\fa,\fb)\in\Sigma_f$, we define $\BD_{\fa,\fb}(a,n)$ to
be the set of all the Brauer $n$-diagrams in $\BD_n(a,n)$ whose
rightmost $f$ horizontal edges in the bottom row are exactly those
pairs in $(\fa,\fb)$. It is clear that $$
\BD_n(a,n)=\bigsqcup_{(\fa,\fb)\in\Sigma_f}\BD_{\fa,\fb}(a,n).
$$
Therefore,
\begin{equation} \label{equa411} v_{\ui}\sum_{D\in\BD_{\fa,\fb}(a,n)}D\neq 0,
\end{equation}
for some $(\fa,\fb)\in\Sigma_f$. Henceforth, we fix such an
$(\fa,\fb)$. We list the integers in the set $\{k|a_1\leq k\leq
n\}\setminus\{a_i, b_i|1\leq i\leq f\}$ as
$t_1,t_2,\cdots,t_{a+1-a_1}$ such that $t_1<t_2<\cdots<t_{a+1-a_1}$.
We define
$$\fI_{\fa,\fb}:=\Bigl\{\uj:=(j_1,\cdots,j_{a+1-a_1})\Bigm|1\leq j_1<\cdots<j_{a+1-a_1}\leq
a\Bigr\}.$$ For each element $\sigma\in\BS_{a+1-a_1}$, we define
$\BD_{\fa,\fb}^{\sigma,\uj}(a,n)$ to be the set of all the Brauer
$n$-diagrams in $\BD_{\fa,\fb}(a,n)$ whose vertical edges connects
the vertex labeled by $t_s$ in the bottom row with the vertex
labeled by $j_{s\sigma}$ in the top row for each integer $s$ with
$1\leq s\leq a+1-a_1$. It is clear that $$
\BD_{(\fa,\fb)}(a,n)=\bigsqcup_{\substack{\sigma\in\BS_{a+1-a_1}\\
\uj\in\fI_{\fa,\fb}}}\BD_{\fa,\fb}^{\sigma,\uj}(a,n).
$$

Now let $\uk:=(k_1,\cdots,k_n)$ be such that $v_{\uk}$ appears
with nonzero coefficient in the expansion of (\ref{equa47}).
Applying Lemma \ref{lm44}, we know that for any $\sigma\in\BS_n$,
\begin{equation}\label{equa410}
\Bigl(\sum_{D\in\BD_n(a,n)}D\Bigr)\sigma=\sum_{D\in\BD_n(a,n)}D.
\end{equation}
As a result, we deduce that $k_1,\cdots,k_n$ are pairwise distinct.
Furthermore, we can assume (if necessary, we replace $v_{\uk}$ by
$v_{\uk}\sigma$ for some $\sigma\in\BS_n$) that
$$
i_{a-s+1}=k_{t_{a+2-a_1-s}},\,\,\text{for $s=1,2,\cdots,a+1-a_1$.}
$$
It follows that $v_{\uk}$ appears with nonzero coefficient in the
expansion of \begin{equation}\label{equa412}
v_{\ui}\sum_{D\in\BD_{\fa,\fb}^{(0)}(a,n)}D,
\end{equation}
where
$\BD_{\fa,\fb}^{(0)}:=\BD_{\fa,\fb}^{1,({a_1},{a_1+1},\cdots,{a})}$.
We divide the remaining proof into two cases:
\smallskip

{\it Case 1.} $\ell_s(\ui)=(n-a)/2$. Then it follows from
(\ref{equa412}) and our results on $\ui$ and $\uk$ that
$\{i_1,\cdots,i_{a_1-1}\}=\{k_1,\cdots,k_{a_1-1}\}$, and the vectors
$v_{i_1},\cdots,v_{i_{a_1-1}}$ lie in a subspace with dimension
$$\leq m-\frac{n-a}{2}-(a+1-a_1)\leq a_1-2.$$ Hence $i_s=i_t$ for some
$1\leq s\neq t\leq a_1-1$. But (\ref{equa412}) implies that
$$ (v_{i_1}\otimes\cdots\otimes
v_{i_{a_1-1}})\sum_{\sigma\in\BS_{a_{1}-1}}\sigma\neq 0,
$$
a contradiction, as required.
\smallskip

{\it Case 2.} $\ell_s(\ui)>(n-a)/2$. Let $W$ be the symplectic
subspace generated by $v_{k_1},\cdots,v_{k_{a_1-1}}$. Note that
$\ell_s(\uk)=\ell_s(\ui)$. It follows that $\dim W\leq 2(a_1-2)$.
Hence $\widetilde{m}:=(\dim W)/2\leq a_1-2$. Note that if
$i_s\not\in W$ for some $1\leq s\leq a_1-1$, then we must have that
$i_s=i'_t$ for some integer $t$ with $1\leq t\leq a_1-1$.
Furthermore, in this case, if we replace the tensor factors
$v_{i_s}, v_{i_t}$ in $v_{\ui}$ by $v_{k_1}, v_{k'_1}$ respectively,
then the coefficient of $v_{\uk}$ in (\ref{equa412}) changes at most
one sign. Therefore, (\ref{equa412}) implies that there exists a
simple $(a_1-1)$-tensor $\widetilde{v}$ in $W^{\otimes a_1-1}$ such
that $v_{k_1}\otimes\cdots\otimes v_{k_{a_1-1}}$ appears with
non-zero coefficient in
$$
\widetilde{v}\sum_{D\in\BD_{a_1-1}(\widetilde{m},\widetilde{m})}D,
$$
where the above element $D$ is understood as element in the Brauer
algebra $\bb_n(-2\widetilde{m})$, acting on the $(a_1-1)$-tensor
space $W^{\otimes a_1-1}$. This is impossible by Proposition
\ref{prop41}. Hence we complete the proof of this
proposition.\end{proof}

\bigskip\bigskip

\bibliographystyle{amsplain}

\end{document}